\documentclass[11pt]{article}
\usepackage{amsfonts}
\usepackage{graphics}
\usepackage{indentfirst}
\usepackage{color}
\usepackage{cite}
\usepackage{latexsym}
\usepackage[paper=a4paper, left=2.1cm, right=2.1cm, top=2.2cm, bottom=1.5cm, headheight=5.5pt, footskip=0.8cm, footnotesep=0.8cm, centering, includefoot]{geometry}
\usepackage{amsmath}
\allowdisplaybreaks
\usepackage{amssymb}
\usepackage[dvips]{epsfig}
\usepackage{amscd}

\newtheorem{theorem}{Theorem}[section]
\newtheorem{remark}{Remark}[section]
\newtheorem{definition}{Definition}[section]
\newtheorem{lemma}{Lemma}[section]

\DeclareMathOperator{\divv}{div}

\makeatletter
\@addtoreset{equation}{section}
\makeatother
\makeatletter
\@addtoreset{equation}{section}
\makeatother

\title{Global strong solution to the two-dimensional density-dependent nematic liquid crystal flows with vacuum}

\author{ Lin Li\thanks{School of Mathematics and Statistics, Chongqing Technology and Business University, Chongqing 400067, People's Republic of China ({\tt lilin420@gmail.com})
}
\quad  Qiao Liu\thanks{Department of Mathematics, Hunan Normal University, Changsha 410081, People's Republic of China ({\tt liuqao2005@163.com}).
}
\quad  Xin Zhong\thanks{Corresponding author. School of Mathematics and Statistics, Southwest University, Chongqing 400715,
People's Republic of China ({\tt xzhong1014@amss.ac.cn}).
}
}
\date{ }

\begin{document}
\maketitle

\begin{abstract}
We are concerned with the Cauchy problem of the two-dimensional (2D) nonhomogeneous incompressible nematic liquid crystal flows on the
whole space $\mathbb{R}^{2}$ with vacuum as far field density. It is proved that the 2D nonhomogeneous incompressible nematic liquid crystal
flows admits a unique global strong solution provided the initial data density and the gradient of orientation decay not too slow at infinity, and the basic energy $\|\sqrt{\rho_0}\mathbf{u}_0\|_{L^2}^2+\|\nabla\mathbf{d}_0\|_{L^2}^2$ is small. In particular, the initial density  may contain vacuum states and even have compact support. Moreover, the large time behavior of the solution is also investigated.
\end{abstract}


Math Subject Classification: 76A15; 35B65; 35Q35

\section{Introduction}
The motion of a two-dimensional nonhomogeneous incompressible nematic liquid crystal flow is governed by the following equations:
\begin{align}
   \label{eq1.1}
 \begin{cases}
\partial_{t}\rho+\divv(\rho\mathbf{u}) =0,\\
\partial_{t}(\rho\mathbf{u})+\divv(\rho\mathbf{u}\otimes\mathbf{u})-\mu\Delta \mathbf{u}+\nabla{P}=-\lambda\divv(\nabla \mathbf{d} \odot\nabla \mathbf{d}),\\
\partial_{t}\mathbf{d}+(\mathbf{u}\cdot\nabla)\mathbf{d}=\gamma(\Delta \mathbf{d}+|\nabla\mathbf{d}|^{2}\mathbf{d}),\\
\quad\divv \mathbf{u}=0,\quad\quad |\mathbf{d}|=1.
\end{cases}
\end{align}
Here, the unknown functions $\rho=\rho(x,t)$, $\mathbf{u}=(u_1,u_2)(x,t)$, and $P=P(x,t)$ denote the density, velocity, and pressure of the fluid, respectively.
$\mathbf{d}=(d_1,d_2,d_3)(x,t)$ is the unknown (averaged) macroscopic/continuum
molecule orientation of the nematic liquid crystal flow. The positive constants $\mu$, $\lambda$ and $\gamma$ represent viscosity, the competition between kinetic energy and potential energy, and $\gamma$ is the microscopic elastic relaxation time for the molecular orientation field, respectively. The notation $\nabla \mathbf{d}\odot\nabla \mathbf{d}$ denotes the $2\times 2$ matrix whose $(i,j)$-th entry is given by $\partial_{i}\mathbf{d}\cdot
\partial_{j}\mathbf{d}$ ($1\leq i,j\leq 2$).

We consider the Cauchy problem for \eqref{eq1.1}  with $(\rho, \mathbf{u})$ vanishing at infinity (in some weak sense). For given initial data $\rho_0$, $\mathbf{u}_0$, and $\mathbf{d}_{0}$, we require that
\begin{align}\label{eq1.2}
 &\rho(x,0)=\rho_{0}(x),\quad \rho\mathbf{u}(x,0)=\rho_0\mathbf{u}_0(x), \quad \mathbf{d}(x,0)=\mathbf{d}_{0}(x),\quad |\mathbf{d}_{0}(x)|=1, \quad\text{ in }\mathbb{R}^{2}.
\end{align}

The above system \eqref{eq1.1}--\eqref{eq1.2} is a macroscopic continuum description of the evolution for
the nematic liquid crystals. It is a simplified
version of the Ericksen-Leslie model \cite{ER,LE}, but it still
retains most important  mathematical structures as well as most of
the essential difficulties of the original Ericksen-Leslie model. A
brief account of the Ericksen-Leslie theory on nematic liquid
crystal flows and the derivations of several approximate systems can
be found in the appendix of \cite{LL1}. For more details on the
hydrodynamic continuum theory of liquid crystals, we refer the
readers to \cite{IWS}. Mathematically, the system \eqref{eq1.1}--\eqref{eq1.2} is a  coupling between the nonhomogeneous incompressible Navier-Stokes equations (see e.g., \cite{ZL,LXZ,Lions1,Temam})
and the transported heat flows of harmonic map (see e.g., \cite{LLZ,W}), and thus, its mathematical analysis is full of challenges.

There is a huge literature on the homogeneous incompressible nematic liquid crystal flows, where namely $\rho$ is constant in \eqref{eq1.1}, refer to \cite{LL1,HMC,LL2,LLW,LW1,LZZ,WD} and references therein. The important progress on the global existence of strong or weak solutions of nonhomogeneous incompressible nematic liquid crystal flows in two dimension has been made recently by some authors. For the initial density away from vacuum, Li \cite{LJK2} established the global existence of strong and weak solutions to the system \eqref{eq1.1}--\eqref{eq1.2} provided that the initial orientation $\mathbf{d}_{0}=(d_{01},d_{02},d_{03})$ satisfies a geometric condition
\begin{equation}\label{lwz}
d_{03}\geq\delta_0\ \text{for some positive}\ \delta_0>0.
\end{equation}
In the presence of vacuum,
if the initial data is small (in some sense) and satisfies the following  compatibility conditions
\begin{align}\label{eq1.3}
-\mu\Delta \mathbf{u}_{0}+\nabla P_{0}+\lambda\divv(\nabla \mathbf{d}_{0}\odot\nabla\mathbf{d}_{0})=\rho_{0}^{\frac{1}{2}} g_{0}
\end{align}
in a bounded smooth domain $\Omega\subseteq\mathbb{R}^2$, and $(P_{0},g_{0})\in H^{1}(\Omega)\times L^{2}(\Omega)$, Wen-Ding \cite{WD} obtained the global existence and uniqueness of the strong solutions to the system \eqref{eq1.1}--\eqref{eq1.2}, see also \cite{FLN,LJK} for related works. It should be emphasized that the possible appearance of vacuum is one
of the main difficulties, which indeed leads to the singular behaviors of solutions in the presence of vacuum, such as the finite time blow-up of smooth solutions \cite{HW1}.

It is not known in general about the existence of global strong solutions to the  problem \eqref{eq1.1}--\eqref{eq1.2} in two-dimension without the geometric condition \eqref{lwz} or the compatibility condition \eqref{eq1.3} imposed on the initial data. This problem is rather interesting and hard to investigate. Indeed, it should be noted that the previous studies on the heat flow of a harmonic map \cite{CDY} indicate that the strong solution of a harmonic map can be blow-up in finite time. In our case, since the system \eqref{eq1.1} contains the heat flow of a harmonic map as a subsystem, we cannot expect that \eqref{eq1.1} have a global strong solution with general initial data. This makes the analysis rather delicate and difficult.

It should be noticed that when $\mathbf{d}$ is a constant vector and $|\mathbf{d}|=1$, the system \eqref{eq1.1} reduces to the nonhomogeneous incompressible Navier-Stokes equations, which have been discussed in
numerous studies \cite{Lions1,LSZ} and so on. It is worth mentioning that L{\"u}-Shi-Zhong \cite{LSZ} recently established the global existence of strong solutions to the 2D Cauchy problem of the incompressible Navier-Stokes equations on the whole space $\mathbb{R}^{2}$ with vacuum as far field density. However, since the system \eqref{eq1.1} contains the incompressible Navier-Stokes equations as a subsystem, one cannot expect, in general, any better results than those for the Navier-Stokes equations. It is a natural and interesting problem to investigate the global existence of strong solutions to the 2D Cauchy problem \eqref{eq1.1}--\eqref{eq1.2} with  vacuum as far field density. In fact, this is the main aim of the present paper.

Before stating the main results, we first explain the notations and
conventions used throughout this paper. For $R>0$, set
\begin{equation*}
B_R \triangleq\left.\left\{x\in\mathbb{R}^2 \right|\,|x|<R \right\},
\quad \int \cdot dx\triangleq\int_{\mathbb{R}^2}\cdot dx.
\end{equation*}
Moreover, for $1\le r\le \infty$ and $k\ge 1$, the standard Lebesgue and  Sobolev spaces are defined as follows:
\begin{equation*}
L^r=L^r(\mathbb{R}^2),\quad W^{k,r}= W^{k,r}(\mathbb{R}^2), \quad H^k = W^{k,2}.
\end{equation*}

Now we define precisely what we mean by strong solutions.\begin{definition}\label{def1.1}
If all derivatives involved in \eqref{eq1.1} for $(\rho, \mathbf{u}, P, \mathbf{d})$ are regular distributions, and equations \eqref{eq1.1} hold almost everywhere in $\mathbb{R}^{2}\times (0,T)$, then $(\rho, \mathbf{u}, P, \mathbf{d})$ is called a strong solution to \eqref{eq1.1}.
\end{definition}

Without loss of generality,  we assume that the initial density $\rho_{0}$ satisfies
\begin{align}\label{eq1.4}
\int_{\mathbb{R}^{2}} \rho_{0} \text{d}x=1,
\end{align}
which implies that there exists a positive constant $N_{0}$ such that
\begin{align}\label{eq1.5}
\int_{B_{N_{0}}}\rho_{0}\text{d}x\geq \frac{1}{2}\int_{\mathbb{R}^{2}}\rho_{0}\text{d}x
= \frac{1}{2},
\end{align}
where $B_{R}\triangleq \{x\in\mathbb{R}^{2} | |x|<R\}$ for all $R>0$. Furthermore, since the concrete values of $\mu$, $\lambda$ and $\gamma$ do
not play a special role in our discussion, in what follows, we assume
\begin{equation*}
\mu=\lambda=\gamma=1
\end{equation*}
throughout this paper.

Now, we state our main result as follows:

\begin{theorem}\label{thm1.2}
For constants $q>2$, $a>1$, assume that the initial data $(\rho_{0}, \mathbf{u}_{0}, \mathbf{d}_{0})$ satisfies \eqref{eq1.4}, \eqref{eq1.5}, and
\begin{align}\label{eq1.6}
\begin{cases}
\rho_{0}\geq 0, \ \rho_{0}\bar{x}^{a}\in L^{1}\cap H^{1}\cap W^{1,q}, \
\sqrt{\rho_{0}}\mathbf{u}_{0}\in L^{2}, \ \nabla \mathbf{u}_{0}\in L^{2},\\
\divv \mathbf{u}_{0}=0,\ \mathbf{d}_{0}\in L^{2},\ \nabla \mathbf{d}_{0} \bar{x}^{\frac{a}{2}}\in L^{2},\ \nabla^{2} \mathbf{d}_{0}\in L^{2},
\ |\mathbf{d}_{0}|=1,
\end{cases}
\end{align}
where
\begin{align*}
\bar{x}\triangleq (e+|x|^{2})^{\frac{1}{2}}\log^{2} (e+|x|^{2}).
\end{align*}
Then there is a positive constant $\varepsilon_0$
depending only on $\|\rho_0\|_{L^1\cap L^\infty},\mu,\lambda,\gamma$
such that if
\begin{equation}\label{z1}
C_{0}\triangleq\int\rho_0|\mathbf{u}_0|^2dx+\int|\nabla\mathbf{d}_0|^2dx
<\varepsilon_0,
\end{equation}
then the Cauchy problem \eqref{eq1.1}--\eqref{eq1.2} has a unique global strong solution $(\rho, \mathbf{u}, P, \mathbf{d})$ satisfying that for any $0<T<\infty$,
\begin{align}\label{eq1.7}
\begin{cases}
0\leq \rho\in C([0,T]; L^{1}\cap H^{1}\cap W^{1,q}),\\
\rho \bar{x}^{a}\in L^{\infty}(0,T; L^{1}\cap H^{1}\cap W^{1,q})\\
\sqrt{\rho}\mathbf{u},  \nabla\mathbf{u}, \sqrt{t}\nabla\mathbf{u},  \sqrt{t}\sqrt{\rho}\mathbf{u},
\sqrt{t}\nabla P, t\nabla P, \sqrt{t}\nabla^{2}\mathbf{u}, t \nabla^{2}\mathbf{u}\in L^{\infty}(0,T; L^{2}),\\
\nabla\mathbf{d}, \nabla\mathbf{d} \bar{x}^{\frac{a}{2}}, \nabla^{2}\mathbf{d}, \sqrt{t}\nabla^{2}\mathbf{d}, \sqrt{t}\nabla \mathbf{d}_{t}, \sqrt{t}\nabla^{3}\mathbf{d},
t\nabla^{3}\mathbf{d}\in L^{\infty}(0,T; L^{2}),\\
\nabla \mathbf{u}\in L^{2}(0,T; H^{1})\cap L^{\frac{q+1}{q}}(0,T; W^{1,q}),\\
\nabla P\in L^{2}(0,T; L^{2})\cap L^{\frac{q+1}{q}}(0,T; L^{q}),\\
\nabla^{2}\mathbf{d}\in L^{2}(0,T; H^{1}),\quad \nabla \mathbf{d}_{t}, \nabla^{2}\mathbf{d} \bar{x}^{\frac{a}{2}}\in L^{2}(\mathbb{R}^{2}\times (0,T)),\\
\sqrt{\rho} \mathbf{u}_{t},  \sqrt{t}\nabla\mathbf{u}_{t}, \sqrt{t}\nabla\mathbf{d}_{t}, \sqrt{t}\nabla^{2}\mathbf{d}_{t}\in L^{2}(\mathbb{R}^{2}\times (0,T)),\\
\sqrt{t} \nabla \mathbf{u} \in L^{2}(0,T; W^{1,q}),
\end{cases}
\end{align}
and
\begin{align}\label{eq1.8}
\inf_{0\leq t\leq T} \int_{B_{N_{1}}} \rho(x,t)\text{d}x\geq \frac{1}{4}
\end{align}
for some positive constant $N_{1}$ depending only on $\|\sqrt{\rho_{0}}\mathbf{u}_{0}\|_{L^{2}}$, $N_{0}$, and $T$. Moreover, the solution
$(\rho, \mathbf{u}, P, \mathbf{d})$ has the following temporal  decay rates, i.e., for all $t\geq 1$,
\begin{align}\label{eq1.9}
\|\nabla \mathbf{u}(\cdot,t)\|_{L^{2}}^2 +\|\nabla^{2} \mathbf{u}(\cdot,t)\|_{L^{2}}+\|\nabla P(\cdot,t)\|_{L^{2}}
+\||\nabla \mathbf{d}||\nabla^{2}\mathbf{d}|(\cdot, t)\|_{L^{2}}+\|\nabla^{2}\mathbf{d}(\cdot, t)\|_{L^{2}}^2
\leq C t^{-1},
\end{align}
where $C$ depends only on $C_0, \|\rho_0\|_{L^1\cap L^\infty}$,  and $\|\nabla \mathbf{u}_0\|_{L^2} $.
\end{theorem}

A few remarks are in order:
\begin{remark}
As stated above, it seems more involved to show the global existence of strong solutions with general initial data. This is the main reason for us to add an additional smallness condition \eqref{z1}. Although it has small energy, its oscillations can be arbitrarily large.
\end{remark}

\begin{remark}
Compared  with \cite{FLN,LJK,WD}, there is no need to impose the additional compatibility condition \eqref{eq1.3} to obtain the global existence of strong solutions.
\end{remark}

\begin{remark}
Our Theorem \ref{thm1.2} holds for the initial density being allowed to have vacuum which is in sharp contrast to \cite{LJK2} where the initial density is absence of vacuum.  Moreover, the geometric condition \eqref{lwz} on the initial orientation is also needed in \cite{LJK2} (see also \cite{LLTZ15}).
\end{remark}

\begin{remark}
It should be noted that our large time decay rates of the velocity and the pressure in \eqref{eq1.9} are the same as those of the incompressible Navier-Stokes equations \cite{LSZ}, hence the orientation has no influence on the large time  behaviors of the velocity and the pressure.
\end{remark}

\begin{remark}
It follows from \eqref{eq1.7} and Aubin-Lions lemma, we see that the velocity $\mathbf{u}$ is continuous with respect to $t$ as long as $t>0$. However, we can not obtain the continuity of the velocity at the initial time due to the presence of vacuum ($\mathbf{u}_0(x)=\mathbf{0}$ if $\rho_0(x)=0$). Nevertheless, we can get the continuity of $\rho\mathbf{u}$ at the initial time.

Indeed, by \eqref{eq1.7}, we immediately have
\begin{align}\label{10}
& \rho\mathbf{u}=\sqrt{\rho}\cdot\sqrt{\rho}\mathbf{u}\in L^\infty_tL^2_x,\\
& \rho\nabla\mathbf{u}\in L^\infty_tL^2_x\label{20},\\
& \nabla\rho\bar{x}^a\in L^\infty_tL^q_x\ \ \text{for}\ q>2\label{30}.
\end{align}
By \eqref{l3} in the next section and \eqref{eq1.7}, we derive that
\begin{equation*}
\mathbf{u}\bar{x}^{-a}\in L^\infty_tL^p_x\ \ \text{for any}\ p>2,
\end{equation*}
which combined with \eqref{30} and H{\"o}lder's inequality leads to
\begin{equation}\label{40}
\nabla\rho\cdot\mathbf{u}=\nabla\rho\bar{x}^a\cdot
\mathbf{u}\bar{x}^{-a}\in L^\infty_tL^2_x.
\end{equation}
Thus, we infer from \eqref{20} and \eqref{40} that
\begin{equation*}
\nabla(\rho\mathbf{u})=\rho\nabla\mathbf{u}+\nabla\rho\cdot\mathbf{u}\in L^\infty_tL^2_x,
\end{equation*}
which along with \eqref{10} yields
\begin{equation}\label{50}
\rho\mathbf{u}\in L^\infty_tH^1_x.
\end{equation}
On the other hand, we deduce from \eqref{l3}, \eqref{eq1.7}, and H{\"o}lder's inequality that
\begin{equation*}
\begin{split}
(\rho\mathbf{u})_t & =\rho_t\mathbf{u}+\rho\mathbf{u}_t \\
& =\nabla\rho|\mathbf{u}|^2+\rho\mathbf{u}_t \\
& =\nabla\rho\bar{x}^{a}\cdot|\mathbf{u}|^2\bar{x}^{-a}
+\sqrt{\rho}\cdot\sqrt{\rho}\mathbf{u}_t \\
& =\nabla\rho\bar{x}^{a}\cdot(|\mathbf{u}|\bar{x}^{-\frac{a}{2}})^2
+\sqrt{\rho}\cdot\sqrt{\rho}\mathbf{u}_t\in L_{t,x}^2,
\end{split}
\end{equation*}
which together with \eqref{50} and Aubin-Lions lemma gives the continuity of $\rho\mathbf{u}$ at the initial time.
\end{remark}

We now make some comments on the analysis of the present paper. Note that for initial data satisfying \eqref{eq1.6}, Liu-Liu-Tan-Zhong \cite{LLTZ15} established the local existence and uniqueness of strong solution to the Cauchy problem \eqref{eq1.1}--\eqref{eq1.2} (see Lemma \ref{lem2.1}). Thus, the proof of Theorem \ref{thm1.2} lies in some global a priori estimates on the strong solutions to the system \eqref{eq1.1}--\eqref{eq1.2} in suitable higher norms.
It should be pointed out that the crucial techniques of proofs in \cite{LW,DHX} cannot be adapted to the situation treated here, since the standard Sobolev embedding inequality is critical in $\mathbb{R}^2$. Moreover, it seems difficult to bound the $L^q(\mathbb{R}^2)$-norm of $\mathbf{u}$ just in terms of $\|\sqrt{\rho}\mathbf{u}\|_{L^{2}(\mathbb{R}^2)}$ and $\|\nabla\mathbf{u} \|_{L^{2}(\mathbb{R}^2)}$.
To this end, we try to adapt some basic ideas used in \cite{LSZ,LXZ}, where the authors investigated the global existence of strong solutions to the Cauchy problem of the 2D nonhomogeneous incompressible Navier--Stokes and MHD equations, respectively. However, compared with \cite{LSZ,LXZ}, for the incompressible nematic liquid crystal flows treated here, the strong coupling terms and strong nonlinear terms, such as $\divv(\nabla\mathbf{d}\odot\nabla\mathbf{d})$,  $\mathbf{u}\cdot\nabla\mathbf{d}$ and $|\nabla\mathbf{d}|^{2}\mathbf{d}$, will bring out some new difficulties.

To overcome these difficulties mentioned above, some new ideas are needed. To deal with the difficulty  caused by the lack of the Sobolev inequality, we observe that, in the momentum equations \eqref{eq1.2}, the velocity $\mathbf{u}$ is always accompanied by $\rho$. Motivated by \cite{LZ,LSZ}, by introducing a weighted function to the density, as well as a Hardy-type inequality (see Lemma \ref{lem2.2}), the $\|\rho^{\eta}\mathbf{u}\|_{L^{r}}$ $(r\geq2,\eta>0)$ is controlled in terms of $\|\sqrt{\rho} \mathbf{u}\|_{L^{2}}$ and $\|\nabla \mathbf{u}\|_{L^{2}}$ (see Lemma \ref{lem2.4}). Then we try to obtain the estimates on the $L^{\infty}(0,T;L^{2})$-norm of $\|\nabla \mathbf{u}\|_{L^{2}}$ and $\|\nabla^{2} \mathbf{d}\|_{L^{2}}$. On the one hand, motivated by \cite{LSZ}, we use material derivatives $\dot{\mathbf{u}}\triangleq \mathbf{u}_{t}+\mathbf{u}\cdot\nabla\mathbf{u}$ instead of the usual $\mathbf{u}_{t}$, and use some facts on Hardy and BMO spaces (see Lemma \ref{lem2.5}) to bound the key term $\int |P||\nabla \mathbf{u}|^{2}\text{d}x$ (see the estimates of $I_{2}$ of \eqref{eq3.6}). On the other hand, the usual $L^{2}(\mathbb{R}^{2}\times (0,T))$-norm of $\nabla \mathbf{d}_{t}$ cannot be directly estimated due to the strong coupled term $\mathbf{u}\cdot\nabla\mathbf{d}$ and the strong nonlinear term $|\nabla \mathbf{d}|^{2}\mathbf{d}$. Motivated by \cite{LXZ1}, multiplying $\eqref{eq3.7}$ by $\Delta\nabla\mathbf{d}$ instead of the usual $\nabla \mathbf{d}_{t}$, and the nonlinear terms $\mathbf{u}\cdot\nabla\mathbf{d}$ and $|\nabla\mathbf{d}|^{2}\mathbf{d}$ can be controlled in terms of $\nabla^{2}\mathbf{d}$ and $\nabla\mathbf{u}$ (see \eqref{eq3.11}), and we find that the key point to obtain the estimate on the $L^\infty(0,T;L^2(\mathbb{R}^2))$-norm of $\nabla\mathbf{u}$ and $\nabla^{2}\mathbf{d}$ is to bound $L^2(0,T;L^2(\mathbb{R}^2))$-norm of $\nabla^{2}\mathbf{d}$ (see \eqref{z3}). Combining the basic energy inequalities (see \eqref{lz} and \eqref{eq3.4}) with Ladyzhenskaya's inequality, we can successfully obtain the a priori bound on the $L^2$-norm of $\nabla^2 \mathbf{d}$ in space and time provided that $\varepsilon_0$ is small (see \eqref{z4} and \eqref{z5}).
Next, using the structure of the 2D heat flows of harmonic maps, we multiply \eqref{eq3.7} by $\nabla\mathbf{d}\Delta|\nabla\mathbf{d}|^{2}$ and thus obtain some useful a priori estimates on $\||\nabla \mathbf{d}||\nabla^{2}\mathbf{d}|\|_{L^{2}}$ and $\||\nabla\mathbf{d}||\Delta\nabla\mathbf{d}|\|_{L^{2}}$ (see \eqref{eq3.17}), which are crucial in deriving the time-independent estimates on both the $L^{\infty}(0,T; L^{2})$-norm of $t^{\frac{1}{2}}\sqrt{\rho}\dot{\mathbf{u}}$ and the $L^{2}(\mathbb{R}^{2}\times (0,T))$-norm of $t^{\frac{1}{2}} \nabla \dot{\mathbf{u}}$ (see \eqref{eq3.13}). By the similar arguments as \cite{LSZ}, we get the bounds of $L^{\infty}(0,T;L^{1})$-norm of spatial weighted estimates of the density (see \eqref{eq3.24}). This together with Lemma \ref{lem2.4} and some careful analysis indicates the desired $L^{1}(0,T;L^{\infty})$ bound of the gradient of the velocity $\nabla\mathbf{u}$ (see \eqref{eq3.31}), which in particular implies the bound on the $L^{\infty}(0,T;L^{q})$-norm ($q>2$) of the gradient of the density. Moreover, some useful spatial weighted estimates on $\rho, \nabla\mathbf{d}, \nabla^{2}\mathbf{d}$ are derived (see Lemma \ref{lem3.6}). With the a priori estimates stated above in hand, we can estimate the higher order derivatives of the solution $(\rho,\mathbf{u},P,\mathbf{d})$ (see \eqref{eq3.37}) by using the same arguments as those in \cite{LSZ,LXZ1} to obtain the desired results.

The remaining parts of the present paper are organized as follows.
In Section 2,  we shall give some elementary facts and inequalities which will be needed in later analysis. In Section 3, we give some a priori estimates which are needed to obtain the global  existence  of strong solutions. Section 4 is devoted to proving Theorem \ref{thm1.2}.

\section{Preliminaries}

In this section, we shall give some known results and elementary inequalities which will be used frequently later.

We start with the local existence of strong solutions whose proof can be found in \cite[Theorem 3.1]{LLTZ15}.

\begin{lemma}\label{lem2.1}
Assume that $(\rho_{0}, \mathbf{u}_{0}, \mathbf{d}_{0})$ satisfies \eqref{eq1.4}--\eqref{eq1.6}. Then there exist a small positive time $T>0$ and a unique strong solution $(\rho, \mathbf{u}, P, \mathbf{d})$ to the Cauchy problem of system \eqref{eq1.1}--\eqref{eq1.2} in $\mathbb{R}^{2}\times (0,T]$ satisfying \eqref{eq1.7} and \eqref{eq1.8}.
\end{lemma}

Next, the following Gagliardo-Nirenberg inequality (see \cite{N1959}) will be used later.
\begin{lemma}[Gagliardo-Nirenberg]\label{lem22}
For $q\in[2,\infty), r\in(2,\infty)$, and $s\in(1,\infty)$, there exists some generic constant $C>0$ which may depend on $q,$ $r$, and $s$ such that for $f\in H^{1}(\mathbb{R}^2)$ and $g\in L^{s}(\mathbb{R}^2)\cap D^{1,r}(\mathbb{R}^2)$, we have
\begin{eqnarray*}
& & \|f\|_{L^q(\mathbb{R}^2)}^{q}\leq C\|f\|_{L^2(\mathbb{R}^2)}^{2}\|\nabla f\|_{L^2(\mathbb{R}^2)}^{q-2}, \\
& & \|g\|_{C(\overline{\mathbb{R}^2})}\leq C\|g\|_{L^s(\mathbb{R}^2)}^{s(r-2)/(2r+s(r-2))}\|\nabla g\|_{L^r(\mathbb{R}^2)}^{2r/(2r+s(r-2))}.
\end{eqnarray*}
\end{lemma}

The following weighted $L^{p}$-bounds for elements of the Hilbert space
$\widetilde{D}^{1,2}(\mathbb{R}^{2})\triangleq \{v\in H^{1}_{loc}(\mathbb{R}^{2}) | \nabla u\in L^{2}(\mathbb{R}^{2})\}$ can be found in Theorem B.1 in \cite{Lions1}.

\begin{lemma}\label{lem2.2}
For $m\in [2,\infty)$ and $\theta\in (\frac{1+m}{2}, \infty)$, there exists a positive constant $C$ such that for  any $v\in \widetilde{D}^{1,2}(\mathbb{R}^{2})$,
\begin{align*}
\left(\int_{\mathbb{R}^{2}}\frac{|v|^{m}}{e+|x|^{2}}(\ln (e+|x|^{2}))^{-\theta}\text{d}x\right)^{\frac{1}{m}}\leq C\|v\|_{L^{2}(B_{1})}+C\|\nabla v\|_{L^{2}(\mathbb{R}^{2})}.
\end{align*}
\end{lemma}

A useful consequence of Lemma \ref{lem2.2} is the following weighted bounds for elements of $\widetilde{D}^{1,2}(\mathbb{R}^{2})$, which have been proved in \cite{LJLZ,LZ,ZL}. It  will play a crucial role in our following analysis.

\begin{lemma}\label{lem2.3}
Let $\bar{x}$ be as in Theorem \ref{thm1.2}. Assume that $\rho\in L^{1}(\mathbb{R}^{2})\cap L^{\infty}(\mathbb{R}^{2})$
is a non-negative function such that
\begin{align*}
\int_{B_{N_{1}}} \rho\text{d}x\geq M_{1}, \quad \|\rho\|_{L^{1}(\mathbb{R}^{2})\cap L^{\infty}(\mathbb{R}^{2})}\leq M_{2},
\end{align*}
for positive constants $M_{1}$, $M_{2}$, and $N_{1}\geq 1$ with $B_{N_{1}}\subseteq \mathbb{R}^{2}$. Then for $\varepsilon >0$ and $\eta>0$, there is a positive constant $C$ depending only on $\varepsilon, \eta, M_{1}, M_{2}, N_{1}$, and $\eta_{0}$ such that every $v\in \widetilde{D}^{1,2}(\mathbb{R}^{2})$ satisfies
\begin{align}\label{l3}
\|v\bar{x}^{-\eta}\|_{L^{\frac{2+\varepsilon}{\widetilde{\eta}}}(\mathbb{R}^{2})}\leq
C\|\sqrt{\rho}  v\|_{L^{2}(\mathbb{R}^{2})}+C\|\nabla v\|_{L^{2}(\mathbb{R}^{2})}
\end{align}
with $\widetilde{\eta}=\min\{1,\eta\}$.
\end{lemma}

\begin{lemma}\label{lem2.4}
Let the assumptions in Lemma \ref{lem2.3} hold. Suppose in addition that $\rho\bar{x}^{a}\in L^{1}(\mathbb{R}^{2})$ with $a>1$. Then for any $\eta\in(0,1]$ and any $s\geq 2$, there is a constant $C$ depending only on $M_{1}, N_{1}, a, \eta, s, \|\rho\|_{L^{\infty}(\mathbb{R}^{2})}$, and $\|\rho\bar{x}^{a}\|_{L^{1}(\mathbb{R}^{2})}$ such that
\begin{align}\label{l2}
\|\rho^{\eta} v\|_{L^{\frac{s}{\eta}}(\mathbb{R}^{2})}\leq C\left(\|\sqrt{\rho} v\|_{L^{2}(\mathbb{R}^{2})}+\|\nabla v\|_{L^{2}(\mathbb{R}^{2})}\right).
\end{align}
\end{lemma}
{\it Proof.} It follows from H{\"o}lder's inequality and Lemma \ref{lem2.3} that
\begin{align*}
\|\rho^\eta v\|_{L^{\frac{s}{\eta}}}
 &  \leq C\|\rho^\eta\bar x^{\frac{3\eta a}{4s}}\|_{L^{\frac{4s}{3\eta}}}
\|v\bar  x^{-\frac{3\eta a}{4s}}\|_{L^{\frac{4s}{\eta}}} \\
 &  \leq C\|\rho\|_{L^\infty}^{\frac{(4s-3)\eta}{4s}}
 \|\rho\bar x^a\|_{L^1}^{\frac{3\eta}{4s}}\left( \|\sqrt{\rho}v\|_{L^2}
 +\|\rho v\|_{L^2}\right) \\
 &  \leq C\left(\|\sqrt{\rho}v\|_{L^2}+\|\nabla v\|_{L^2}\right),
\end{align*}
which implies \eqref{l2}.  \hfill $\Box$

Finally, let $\mathcal{H}^{1}(\mathbb{R}^2)$ and BMO$(\mathbb{R}^2)$ stand for the usual Hardy and BMO spaces (see \cite[Chapter IV]{Stein}). Then the following well-known facts play a key role in the proof of Lemma \ref{lem3.2} in the next section.
\begin{lemma}\label{lem2.5}
(a) There is a positive constant $C$ such that
\begin{equation*}
\|\mathbf{E}\cdot \mathbf{B}\|_{\mathcal{H}^{1}}
\leq C\|\mathbf{E}\|_{L^{2}}\|\mathbf{B}\|_{L^{2}}
\end{equation*} for all $\mathbf{E}\in L^{2}(\mathbb{R}^2)$ and $\mathbf{B}\in L^{2}(\mathbb{R}^2)$ satisfying
\begin{equation*}
\divv \mathbf{E}=0,\ \nabla^{\bot}\cdot \mathbf{B}=0\ \ \text{in}\ \ \mathcal{D}'(\mathbb{R}^2).
\end{equation*}
(b) There is a positive constant $C$ such that
\begin{equation}\label{lem1}
\| v \|_{{\rm BMO}}\leq C\|\nabla v \|_{L^{2}}
\end{equation} for all $  v \in \tilde{D}^{1,2}(\mathbb{R}^2)$.
\end{lemma}
{\it Proof.}
(a) For the detailed proof, see \cite[Theorem II.1]{CLMS93}.

(b) It follows  from the Poincar{\'e} inequality that for any ball $B\subset\mathbb{R}^2$
\begin{equation*}
\frac{1}{|B|}\int_{B}\left| v(x) - \frac{1}{|B|}\int_B v(y)dy\right|dx\leq C\left(\int_{B}|\nabla v |^2dx\right)^{1/2},
\end{equation*}
which directly gives \eqref{lem1}.  \hfill $\Box$

\section{A priori estimates}\label{sec3}
In this section, we shall establish some necessary a priori estimates for strong solutions $(\rho, \mathbf{u}, P, \mathbf{d})$ to the Cauchy problem of system \eqref{eq1.1}--\eqref{eq1.2} to extend the local strong  solutions guaranteed by Lemma \ref{lem2.1}.  Thus, let $T>0$ be a fixed time and $(\rho, \mathbf{u}, P, \mathbf{d})$ be the strong solution to system \eqref{eq1.1}--\eqref{eq1.2} on $\mathbb{R}^{2}\times (0,T]$ with initial data $(\rho_{0}, \mathbf{u}_{0}, \mathbf{d}_{0})$ satisfying \eqref{eq1.4}--\eqref{z1}. In what follows, the convention of summation over repeated indices is used.

\subsection{Lower order estimates}

\begin{lemma}\label{lem3.2}
If $\varepsilon_0$ in \eqref{z1} depending only on $\|\rho_0\|_{L^1\cap L^\infty},\mu,\lambda,\gamma$ is sufficiently small, then there exists a positive constant $C$  depending only on $C_0$, $\|\rho_{0}\|_{L^{1}\cap L^{\infty}}$,
$\|\nabla \mathbf{u}_{0}\|_{L^{2}}$, and $\|\nabla^{2} \mathbf{d}_{0}\|_{L^{2}}$ such that
\begin{align}\label{eq3.2}
&\sup_{0\leq t\leq T}\left(\|\rho\|_{L^{1}\cap L^{\infty}}+\|\nabla \mathbf{u}\|_{L^{2}}^{2} +\|\nabla^{2} \mathbf{d}\|_{L^{2}}^{2}\right)
+\int_{0}^{T} \left( \|\sqrt{\rho} \dot{\mathbf{u}}\|_{L^{2}}^{2} +\|\nabla^{3} \mathbf{d}\|_{L^{2}}^{2}\right)\text{d}t
\leq C,
\end{align}
where $\dot{\mathbf{u}}\triangleq \mathbf{u}_{t}+\mathbf{u}\cdot\nabla \mathbf{u}$. Furthermore,  we have
\begin{align}\label{eq3.3}
&\sup_{0\leq t\leq T} t\left( \|\nabla \mathbf{u}\|_{L^{2}}^{2} +\|\nabla^{2} \mathbf{d}\|_{L^{2}}^{2}\right)
+\int_{0}^{T} t\left( \|\sqrt{\rho} \dot{\mathbf{u}}\|_{L^{2}}^{2} +\|\nabla^{3} \mathbf{d}\|_{L^{2}}^{2}\right)\text{d}t
\leq C.
\end{align}
\end{lemma}

{\it Proof.}
First, since $\divv \mathbf{u}=0$, it is easy to obtain from the equation $\eqref{eq1.1}_{1}$ (see \cite{Lions1}) that
\begin{align}\label{eq3.5}
\|\rho(t)\|_{L^{p}}=\|\rho_{0}\|_{L^{p}} \ \text{ for all } p\in [1,\infty] \text{ and } t\geq 0.
\end{align}
Applying standard energy estimate to \eqref{eq1.1} (see \cite[lemma 3.1]{LJK}) gives that
\begin{align}\label{lz}
&\sup_{0\leq t\leq T} \left(
\|\sqrt{\rho} \mathbf{u}\|_{L^{2}}^{2}+\|\nabla \mathbf{d}\|_{L^{2}}^{2}\right)
+\int_{0}^{T}\|\nabla \mathbf{u}\|_{L^{2}}^{2}\text{d}t
\leq C_0,
\end{align}
and
\begin{align}\label{eq3.4}
\sup_{0\leq t\leq T} \left(
\|\sqrt{\rho} \mathbf{u}\|_{L^{2}}^{2}+\|\nabla \mathbf{d}\|_{L^{2}}^{2}\right)
+\int_{0}^{T}\left(\|\nabla \mathbf{u}\|_{L^{2}}^{2}+\|\Delta \mathbf{d}\|_{L^{2}}^{2}\right)\text{d}t
\leq C_0+C\int_{0}^{T} \|\nabla \mathbf{d}\|_{L^{4}}^{4}\text{d}t.
\end{align}
Now, motivated by \cite{LSZ},  multiplying  $\eqref{eq1.1}_{2}$ by $\dot{\mathbf{u}}$ and then integrating the resulting equality over $\mathbb{R}^2$ lead to
\begin{align}\label{eq3.6}
\int \rho |\dot{\mathbf{u}}|^{2}\text{d}x=&\int \Delta \mathbf{u}\cdot\dot{\mathbf{u}}\text{d}x -\int \nabla P\cdot\dot{\mathbf{u}}\text{d}x-\int \divv(\nabla \mathbf{d}\odot\nabla \mathbf{d}) \cdot\dot{\mathbf{u}}\text{d}x\nonumber\\
\triangleq& I_{1}+I_{2}+I_{3}.
\end{align}
It follows from integration by parts and Gagliardo-Nirenberg inequality that
\begin{align*}
I_{1}  &  =\int\Delta \mathbf{u} \cdot(\partial_{t} \mathbf{u} + \mathbf{u} \cdot\nabla \mathbf{u} )dx \notag \\
 &  = -\frac{1}{2}\frac{d}{dt}\|\nabla \mathbf{u} \|_{L^2}^{2}
-\int\partial_{i}u_{j}\partial_{i}(u_{k}\partial_{k}u_{j})dx \notag \\
 &  \leq-\frac{1}{2}\frac{d}{dt}\|\nabla \mathbf{u} \|_{L^2}^{2}+C\|\nabla\mathbf{u} \|_{L^3}^{3} \notag \\
 &  \leq-\frac{1}{2}\frac{d}{dt}\|\nabla\mathbf{u}\|_{L^2}^{2}
+C\|\nabla\mathbf{u}\|_{L^2}^{2}\|\nabla^{2}\mathbf{u}\|_{L^2}.
\end{align*}
We deduce from integration by parts and \eqref{eq1.1}$_4$ that
\begin{align*}
I_{2}  &  =-\int\nabla P\cdot(\partial_{t}\mathbf{u}+ \mathbf{u}\cdot\nabla\mathbf{u})dx \notag \\
 &  = \int P\partial_{j}u_{i}\partial_{i}u_{j}dx \notag \\
  & \leq C\|P\|_{\rm BMO}\|\partial_{j}u_{i}\partial_{i}u_{j}\|_{\mathcal{H}^1},
\end{align*}
where one has used the duality of $\mathcal{H}^1$ space and  BMO  one (see \cite[Charpter IV]{Stein}) in the last inequality. Since $\divv(\partial_j \mathbf{u} )=\partial_j\divv \mathbf{u} =0$ and $\nabla^{\bot}\cdot(\nabla u_{j})=0$,   Lemma \ref{lem2.5} yields
\begin{equation*}
|I_{2}|=\left|\int P\partial_{j}u_{i}\partial_{i}u_{j}dx \right|\leq C\|\nabla P\|_{L^2}\|\nabla\mathbf{u}\|_{L^2}^{2}.
\end{equation*}
To bound the term $I_{3}$, we first apply  $\nabla $ on $\eqref{eq1.1}_{3}$ to get
\begin{equation}\label{eq3.7}
\nabla \mathbf{d}_{t}-\Delta\nabla \mathbf{d}=-\nabla(\mathbf{u}\cdot \nabla \mathbf{d})+\nabla(|\nabla \mathbf{d}|^{2}\mathbf{d}),
\end{equation}
which combined with  H\"{o}lder's and Gagliardo-Nirenberg inequalities leads to
\begin{align*}
I_{3}= & \int (\nabla \mathbf{d}\odot\nabla \mathbf{d})\cdot\nabla \mathbf{u}_{t}\text{d}x -\int \divv(\nabla\mathbf{d}\odot\nabla \mathbf{d})\cdot(\mathbf{u}\cdot\nabla \mathbf{u})\text{d}x\nonumber\\
= & \frac{d}{dt}\!\int \!(\nabla \mathbf{d}\odot\nabla \mathbf{d})\cdot\nabla \mathbf{u}\text{d}x -\!\!\int\! (\nabla \mathbf{d}_{t}\odot\nabla \mathbf{d})\cdot\nabla \mathbf{u}\text{d}x -\!\!\int\! (\nabla \mathbf{d}\odot\nabla \mathbf{d}_{t})\cdot\nabla \mathbf{u}\text{d}x\nonumber\\
& -\!\!\int\!\divv(\nabla \mathbf{d}\odot\nabla \mathbf{d})\cdot(\mathbf{u}\cdot\nabla \mathbf{u})\text{d}x\nonumber\\
= & \frac{d}{dt}\!\int \!(\nabla \mathbf{d}\odot\nabla \mathbf{d})\cdot\nabla \mathbf{u}\text{d}x -\!\!\int\! [(\Delta\nabla \mathbf{d}-\nabla(\mathbf{u}\cdot \nabla \mathbf{d})+\nabla(|\nabla \mathbf{d}|^{2}\mathbf{d}))\odot\nabla \mathbf{d}]\cdot\nabla \mathbf{u}\text{d}x\nonumber\\
& -\int\! [\nabla \mathbf{d}\odot(\Delta\nabla \mathbf{d}-\nabla(\mathbf{u}\cdot \nabla \mathbf{d})+\nabla(|\nabla \mathbf{d}|^{2}\mathbf{d}))]\cdot\nabla \mathbf{u}\text{d}x
-\!\!\int\!\divv(\nabla \mathbf{d}\odot\nabla \mathbf{d})\cdot(\mathbf{u}\cdot\nabla\mathbf{u})\text{d}x\nonumber\\
= & \frac{d}{dt}\!\int \!(\nabla \mathbf{d}\odot\!\nabla \mathbf{d})\cdot\!\nabla \mathbf{u}\text{d}x -\!\!\int\! [(\Delta\nabla \mathbf{d}-\!\nabla \mathbf{u}\cdot \nabla \mathbf{d}+\!\nabla(|\nabla \mathbf{d}|^{2}\mathbf{d}))\odot\!\nabla \mathbf{d}]\cdot\nabla \mathbf{u}\text{d}x\nonumber\\
& + \int\! u_{k}\partial_{k}\partial_{i}d_{\ell} \partial_{j}d_{\ell}\partial_{j}u_{i}\text{d}x
-\int\! [\nabla \mathbf{d}\odot(\Delta\nabla \mathbf{d}-\nabla\mathbf{u}\cdot \nabla \mathbf{d}+\nabla(|\nabla \mathbf{d}|^{2}\mathbf{d}))]\cdot\nabla \mathbf{u}\text{d}x\nonumber\\
& + \int \partial_{i}d_{\ell} u_{k}\partial_{k} \partial_{j}d_{\ell} \partial_{j} u_{i}\text{d}x+\int\!\partial_{i} d_{\ell}\partial_{j} d_{\ell}\partial_{j}(u_{k}\partial_{k} u_{i})\text{d}x\nonumber\\
= & \frac{d}{dt}\!\int \!(\nabla \mathbf{d}\odot\nabla \mathbf{d})\cdot\nabla \mathbf{u}\text{d}x - \!\!\int\! [(\Delta\nabla \mathbf{d}-\!\nabla \mathbf{u}\cdot \nabla \mathbf{d}+\!\nabla(|\nabla \mathbf{d}|^{2}\mathbf{d}))\odot\!\nabla \mathbf{d}]\cdot\nabla \mathbf{u}\text{d}x\nonumber\\
&-\int\! [\nabla \mathbf{d}\odot(\Delta\nabla \mathbf{d}-\nabla\mathbf{u}\cdot \nabla \mathbf{d}+\nabla(|\nabla \mathbf{d}|^{2}\mathbf{d}))]\cdot\nabla \mathbf{u}\text{d}x
+\int \partial_{i}d_{\ell}\partial_{j}d_{\ell} \partial_{j} u_{k}\partial_{k} u_{i}\text{d}x\nonumber\\
\leq & \frac{d}{dt}\!\int \!(\nabla \mathbf{d}\odot\nabla \mathbf{d})\cdot\nabla \mathbf{u}\text{d}x +C\|\nabla^{3} \mathbf{d}\|_{L^{2}} \|\nabla \mathbf{d}\|_{L^{6}} \|\nabla \mathbf{u}\|_{L^{3}}
+C\|\nabla \mathbf{d}\|_{L^{6}}^{2}\|\nabla \mathbf{u}\|_{L^{3}}^{2}\nonumber\\
&+C\|\nabla \mathbf{d}\|_{L^{6}}^{4} \|\nabla \mathbf{u}\|_{L^{3}} +C\|\nabla^{2}\mathbf{d}\|_{L^{3}}\|\nabla \mathbf{d}\|_{L^{6}}^{2} \|\nabla \mathbf{u}\|_{L^{3}}\nonumber\\
\leq& \frac{d}{dt}\!\int \!(\nabla \mathbf{d}\odot\nabla \mathbf{d})\cdot\nabla \mathbf{u}\text{d}x +\frac{\varepsilon}{4}\|\nabla^{3} \mathbf{d}\|_{L^{2}}^{2}+C\|\nabla \mathbf{u}\|_{L^{3}}^{3}+C\|\nabla \mathbf{d}\|_{L^{6}}^{6}+C\|\nabla^{2} \mathbf{d}\|_{L^{3}}^{\frac{3}{2}} \|\nabla \mathbf{d}\|_{L^{6}}^{3}\ (\text{with } \varepsilon>0)\nonumber\\
\leq & \frac{d}{dt}\!\int \!(\nabla \mathbf{d}\odot\nabla \mathbf{d})\cdot\nabla \mathbf{u}\text{d}x +\frac{\varepsilon}{4}\|\nabla^{3} \mathbf{d}\|_{L^{2}}^{2}+C\|\nabla \mathbf{u}\|_{L^{3}}^{3}+C\|\nabla \mathbf{d}\|_{L^{6}}^{6}+C\|\nabla^{3} \mathbf{d}\|_{L^{2}}^{\frac{3}{4}} \|\nabla \mathbf{d}\|_{L^{6}}^{\frac{15}{4}}\quad\nonumber\\
\leq & \frac{d}{dt}\!\int \!(\nabla \mathbf{d}\odot\nabla \mathbf{d})\cdot\nabla \mathbf{u}\text{d}x +\frac{\varepsilon}{2}\|\nabla^{3} \mathbf{d}\|_{L^{2}}^{2}+C\|\nabla \mathbf{u}\|_{L^{3}}^{3}+C\|\nabla \mathbf{d}\|_{L^{6}}^{6}\nonumber\\
\leq & \frac{d}{dt}\!\int \!(\nabla \mathbf{d}\odot\nabla \mathbf{d})\cdot\nabla \mathbf{u}\text{d}x +\frac{\varepsilon}{2}\|\nabla^{3} \mathbf{d}\|_{L^{2}}^{2}+C\|\nabla \mathbf{u}\|_{L^{2}}^{2}
\|\nabla^{2} \mathbf{u}\|_{L^{2}}+C\|\nabla \mathbf{d}\|_{L^{2}}^{2}\|\nabla^{2} \mathbf{d}\|_{L^{2}}^{4}\\
\leq & \frac{d}{dt}\!\int \!(\nabla \mathbf{d}\odot\nabla \mathbf{d})\cdot\nabla \mathbf{u}\text{d}x +\frac{\varepsilon}{2}\|\nabla^{3} \mathbf{d}\|_{L^{2}}^{2}+C\|\nabla \mathbf{u}\|_{L^{2}}^{2}
\|\nabla^{2} \mathbf{u}\|_{L^{2}}+CC_{0}\|\nabla^{2} \mathbf{d}\|_{L^{2}}^{4},
\end{align*}
where in the last inequality we have used \eqref{lz}.
Inserting the above estimates of $I_{i}\ (i=1,2,3)$ into \eqref{eq3.6}, and then using \eqref{eq3.4},  it holds that
\begin{align}\label{eq3.8}
&\frac{d}{dt}\left(\frac{1}{2}\|\nabla \mathbf{u}\|_{L^{2}}^{2}-\int(\nabla d\odot\nabla \mathbf{d})\cdot\nabla \mathbf{u}\text{d}x\right)+\|\sqrt{\rho} \dot{\mathbf{u}}\|_{L^{2}}^{2}\nonumber\\
& \leq \frac{\varepsilon}{2}\|\nabla^{3} \mathbf{d}\|_{L^{2}}^{2}+CC_{0}\|\nabla^{2} \mathbf{d}\|_{L^{2}}^{4}
+C(\|\nabla^{2} \mathbf{u}\|_{L^{2}}+\|\nabla P\|_{L^{2}})\|\nabla \mathbf{u}\|_{L^{2}}^{2}.
\end{align}
On the other hand, since $(\rho, \mathbf{u}, P, \mathbf{d})$ satisfies the following Stokes system
\begin{align*}
\begin{cases}
-\Delta\mathbf{u}+\nabla P =-\rho\dot{\mathbf{u}} -\divv(\nabla \mathbf{d}\odot\nabla \mathbf{d}),\quad& x\in\mathbb{R}^{2},\\
\divv\mathbf{u}=0,\quad& x\in\mathbb{R}^{2},\\
\mathbf{u}(x)\rightarrow \mathbf{0},\quad& |x|\rightarrow \infty,
\end{cases}
\end{align*}
applying the standard $L^{p}$-estimate to the above system (see \cite{Temam}) gives that for any $p\in(1,\infty)$,
\begin{align}\label{eq3.9}
\|\nabla^{2}\mathbf{u}\|_{L^{p}}+\|\nabla P\|_{L^{p}}\leq C(p)\left(\|\rho \dot{\mathbf{u}}\|_{L^{p}}+\||\nabla \mathbf{d}||\Delta \mathbf{d}|\|_{L^{p}}\right)
\leq C(p)\left(\|\sqrt{\rho} \dot{\mathbf{u}}\|_{L^{p}}+\||\nabla \mathbf{d}||\nabla^{2} \mathbf{d}|\|_{L^{p}}\right),
\end{align}
where we have used the identity $\divv(\nabla\mathbf{d}\odot\nabla\mathbf{d})
=\nabla\mathbf{d}\cdot\Delta\mathbf{d}$ and \eqref{eq3.5}. In particular, we derive
\begin{align}\label{l1}
\|\nabla^{2}\mathbf{u}\|_{L^{2}}+\|\nabla P\|_{L^{2}}
& \leq C\left(\|\sqrt{\rho} \dot{\mathbf{u}}\|_{L^{2}}+\||\nabla \mathbf{d}||\nabla^{2} \mathbf{d}|\|_{L^{2}}\right)\nonumber\\
& \leq C\left(\|\sqrt{\rho} \dot{\mathbf{u}}\|_{L^{2}}+\|\nabla \mathbf{d}\|_{L^{4}}\|\nabla^{2}\mathbf{d}\|_{L^{4}}\right)\nonumber\\
& \leq C\left(\|\sqrt{\rho}\dot{\mathbf{u}}\|_{L^{2}}+\|\nabla \mathbf{d}\|_{L^{2}}^{\frac{1}{2}} \|\nabla^{2}\mathbf{d}\|_{L^{2}}\|\nabla^{3}\mathbf{d}\|_{L^{2}}^{\frac{1}{2}}
\right),
\end{align}
which combined with \eqref{eq3.8} and Young's inequality leads to
\begin{align}\label{eq3.10}
\frac{d}{dt} B(t) +\|\sqrt{\rho}\dot{\mathbf{u}}\|_{L^{2}}^{2}
\leq \varepsilon \|\nabla^{3} \mathbf{d}\|_{L^{2}}^{2}+\varepsilon \|\sqrt{\rho} \dot{\mathbf{u}}\|_{L^{2}}^{2}
+C\|\nabla\mathbf{u}\|_{L^{2}}^{4}+CC_{0}\|\nabla^{2}\mathbf{d}\|_{L^{2}}^{4},
\end{align}
where
\begin{align*}
B(t)\triangleq\frac{1}{2} \|\nabla\mathbf{u}\|_{L^{2}}^{2}-\int (\nabla \mathbf{d}\odot\nabla \mathbf{d})\cdot\nabla \mathbf{u}\text{d}x
\end{align*}
satisfies
\begin{align}\label{z2}
\frac{1}{4} \|\nabla\mathbf{u}\|_{L^{2}}^{2} -C_{1}C_{0} \|\nabla^{2} \mathbf{d}\|_{L^{2}}^{2}\leq B(t)\leq C\|\nabla \mathbf{u}\|_{L^{2}}^{2}+CC_0\|\nabla^{2} \mathbf{d}\|_{L^{2}}^{2}
\end{align}
owing to the following estimate
\begin{align*}
\left|\int (\nabla \mathbf{d}\odot\nabla \mathbf{d})\cdot\nabla\mathbf{u}\text{d}x\right|
\leq & \frac{1}{4} \|\nabla\mathbf{u}\|_{L^{2}}^{2}+C\|\nabla \mathbf{d}\|_{L^{4}}^{4} \nonumber\\
\leq & \frac{1}{4} \|\nabla \mathbf{u}\|_{L^{2}}^{2}+C\|\nabla \mathbf{d}\|_{L^{2}}^{2} \|\nabla^{2}\mathbf{d}\|_{L^{2}}^{2}\nonumber\\
\leq & \frac{1}{4} \|\nabla \mathbf{u}\|_{L^{2}}^{2}+CC_{0} \|\nabla^{2}\mathbf{d}\|_{L^{2}}^{2}.
\end{align*}

Now, multiplying $\eqref{eq3.7}$ by $-\nabla\Delta \mathbf{d}$ and then integrating by parts over $\mathbb{R}^{2}$, it follows from H\"{o}lder's and Gagliardo-Nirenberg inequalities, \eqref{lz}, and \eqref{eq3.9} that
\begin{align}\label{eq3.11}
& \frac{1}{2}\frac{d}{dt}\|\nabla^{2} \mathbf{d}\|_{L^{2}}^{2}+\|\nabla\Delta\mathbf{d}\|_{L^{2}}^{2}\nonumber\\
& =  \int \nabla(\mathbf{u}\cdot\nabla \mathbf{d})\cdot\nabla\Delta \mathbf{d}\text{d}x
-\int \nabla (|\nabla \mathbf{d}|^{2} \mathbf{d})\nabla\Delta \mathbf{d}\text{d}x\nonumber\\
& = \int (\nabla \mathbf{u}\cdot\nabla \mathbf{d})\cdot \nabla\Delta \mathbf{d}\text{d}x+\int u_{i}\partial_{i}\partial_{j} d_{\ell}\partial_{i}\partial_{kk}d_{\ell}\text{d}x- \int \nabla (|\nabla \mathbf{d}|^{2} \mathbf{d})\nabla\Delta \mathbf{d}\text{d}x\nonumber\\
& = \int (\nabla \mathbf{u}\cdot\nabla \mathbf{d})\cdot \nabla\Delta \mathbf{d}\text{d}x-\int \partial_{k} u_{i} \partial_{i}\partial_{j} d_{\ell}
\partial_{i} \partial_{k} d_{\ell} \text{d}x-\int \nabla(|\nabla \mathbf{d}|^{2}\mathbf{d})\cdot\nabla\Delta \mathbf{d}\text{d}x\nonumber\\
& \leq  C\!\left(\|\nabla^{3} \mathbf{d}\|_{\!L^{2}} \|\nabla \mathbf{u}\|_{\!L^{3}} \|\nabla \mathbf{d}\|_{\!L^{6}}\! +\|\nabla \mathbf{u}\|_{\!L^{3}}\|\nabla^{2} \mathbf{d}\|_{\!L^{3}}^{2}\!
+
\|\nabla^{3} \mathbf{d}\|_{\!L^{2}} \|\nabla \mathbf{d}\|_{\!L^{6}}^{3}\! +\|\nabla^{3} \mathbf{d}\|_{\!L^{2}}\|\nabla^{2}\mathbf{d}\|_{\!L^{3}} \|\nabla \mathbf{d}\|_{\!L^{6}}\right)\nonumber\\
& \leq \frac{\varepsilon}{4} \|\nabla^{3} \mathbf{d}\|_{L^{2}}^{2} +C\|\nabla \mathbf{u}\|_{L^{3}}^{3}+C\|\nabla \mathbf{d}\|_{L^{6}}^{6}+C\|\nabla^{2}\mathbf{d}\|_{L^{3}}^{3}+C\|\nabla^{2} \mathbf{d}\|_{L^{3}}^{2}\|\nabla \mathbf{d}\|_{L^{6}}^{2}\nonumber\\
& \leq  \frac{\varepsilon}{4} \|\nabla^{3} \mathbf{d}\|_{L^{2}}^{2} +C\|\nabla \mathbf{u}\|_{L^{3}}^{3}+C\|\nabla \mathbf{d}\|_{L^{6}}^{6}+C\|\nabla^{3}\mathbf{d}\|_{L^{2}}^{\frac{3}{2}}\|\nabla \mathbf{d}\|_{L^{2}}^{\frac{3}{2}}+C\|\nabla^{3}\mathbf{d}\|_{L^{2}}\|\nabla \mathbf{d}\|_{L^{6}}^{3}\nonumber\\
& \leq  \frac{\varepsilon}{2} \|\nabla^{3} \mathbf{d}\|_{L^{2}}^{2} +C\|\nabla \mathbf{u}\|_{L^{3}}^{3}+C\|\nabla \mathbf{d}\|_{L^{6}}^{6}\nonumber\\
& \leq  \frac{\varepsilon}{2} \|\nabla^{3} \mathbf{d}\|_{L^{2}}^{2} +C\|\nabla \mathbf{u}\|_{L^{2}}^{2}\|\nabla^{2} \mathbf{u}\|_{L^{2}}
+C \|\nabla \mathbf{d}\|_{L^{2}}^{2}\|\nabla^{2}\mathbf{d} \|_{L^{2}}^{4}\nonumber\\
& \leq  \frac{\varepsilon}{2} \|\nabla^{3} \mathbf{d}\|_{L^{2}}^{2}+\frac{\varepsilon}{2} (\|\sqrt{\rho} \dot{\mathbf{u}}\|_{L^{2}}^{2}+\||\nabla \mathbf{d}||\nabla^{2}\mathbf{d}|\|_{L^{2}}^{2})+CC_{0}\|\nabla^{2} \mathbf{d}\|_{L^{2}}^{4} +
C\|\nabla \mathbf{u}\|_{L^{2}}^{4}\nonumber\\
& \leq \varepsilon\|\nabla^{3} \mathbf{d}\|_{L^{2}}^{2}+\varepsilon \|\sqrt{\rho} \dot{\mathbf{u}}\|_{L^{2}}^{2}
+CC_{0}\|\nabla^{2} \mathbf{d}\|_{L^{2}}^{4} +
C\|\nabla \mathbf{u}\|_{L^{2}}^{4}.
\end{align}
Multiplying \eqref{eq3.11} by $2(C_{1}C_{0}+1)$, then adding the resulting inequality with \eqref{eq3.10} and choosing $\varepsilon$ suitably small, we obtain
\begin{align}\label{eq3.12}
& \frac{d}{dt}\left(B(t)+(C_{1}C_{0}+1)\|\nabla^{2} \mathbf{d}\|_{L^{2}}^{2}\right)+\frac{1}{2}\|\sqrt{\rho} \dot{\mathbf{u}}\|_{L^{2}}^{2}+\|\nabla^{3}\mathbf{d}\|_{L^{2}}^{2}\nonumber\\
& \leq C\left(\|\nabla^{2} \mathbf{d}\|_{L^{2}}^{2}\|\nabla^{2} \mathbf{d}\|_{L^{2}}^{2} +
\|\nabla \mathbf{u}\|_{L^{2}}^{2}\|\nabla \mathbf{u}\|_{L^{2}}^{2}\right).
\end{align}
This along with \eqref{lz}, \eqref{z2}, and Gronwall's inequality yields
\begin{align}\label{z3}
&\sup_{0\leq t\leq T}\left(\|\nabla \mathbf{u}\|_{L^{2}}^{2} +\|\nabla^{2} \mathbf{d}\|_{L^{2}}^{2}\right)
+\int_{0}^{T} \left( \|\sqrt{\rho} \dot{\mathbf{u}}\|_{L^{2}}^{2} +\|\nabla^{3}\mathbf{d}\|_{L^{2}}^{2}\right)\text{d}t
\leq C+C\exp \left[C\int_{0}^{T}\|\nabla^{2}\mathbf{d}\|_{L^{2}}^{2}\text{d}t\right].
\end{align}
It follows from \eqref{eq3.4}, \eqref{lz}, and Ladyzhenskaya's inequality that
\begin{align}\label{z4}
\int_{0}^{T}\|\nabla^{2}\mathbf{d}\|_{L^{2}}^{2}\text{d}t
&\leq CC_{0}+C\int_{0}^{T}\|\nabla\mathbf{d}\|_{L^{4}}^{4}\text{d}t\nonumber\\
&\leq CC_{0}+C\sup_{0\leq t\leq T}\|\nabla\mathbf{d}\|_{L^{2}}^{2}\int_{0}^{T}
\|\nabla^2\mathbf{d}\|_{L^{2}}^{2}\text{d}t\nonumber\\
&\leq CC_{0}+CC_{0}\int_{0}^{T}
\|\nabla^2\mathbf{d}\|_{L^{2}}^{2}\text{d}t,
\end{align}
and thus
\begin{align}\label{z5}
\int_{0}^{T}\|\nabla^{2}\mathbf{d}\|_{L^{2}}^{2}\text{d}t
\leq C
\end{align}
provided $\varepsilon_0$ in \eqref{z1} is small. In particular,
\begin{equation*}
\int_{0}^{T}\|\nabla\mathbf{d}\|_{L^{4}}^{4}\text{d}t
\leq CC_{0}+CC_{0}\int_{0}^{T}
\|\nabla^2\mathbf{d}\|_{L^{2}}^{2}\text{d}t
\leq C,
\end{equation*}
which implies
\begin{align}\label{ll}
\sup_{0\leq t\leq T} \left(
\|\sqrt{\rho} \mathbf{u}\|_{L^{2}}^{2}+\|\nabla \mathbf{d}\|_{L^{2}}^{2}\right)
+\int_{0}^{T}\left(\|\nabla \mathbf{u}\|_{L^{2}}^{2}+\|\Delta \mathbf{d}\|_{L^{2}}^{2}\right)\text{d}t
\leq C.
\end{align}
Combining \eqref{z3} with \eqref{z5}, we derive \eqref{eq3.2}.

Finally,  multiplying \eqref{eq3.12} by $t$, and then applying Gronwall's inequality to the resulting inequality, it follows from \eqref{z5} and \eqref{z2} gives \eqref{eq3.3}. This completes the proof of Lemma \ref{lem3.2}.  \hfill $\Box$

Next, motivated by \cite{LSZ,LXZ1}, we have the following estimates on the material derivatives of the velocity which are important for the higher order estimates of both the density and the velocity.
\begin{lemma}\label{lem3.3}
There exists a positive constant $C$ depending only on $C_0$, $\|\rho_{0}\|_{L^{1}\cap L^{\infty}}$,
$\|\nabla \mathbf{u}_{0}\|_{L^{2}}$, and $\|\nabla^{2} \mathbf{d}_{0}\|_{L^{2}}$ such that for $i=1,2$
\begin{equation}\label{eq3.13}
\sup_{0\leq t\leq T}t^{i}\left(\|\sqrt{\rho} \dot{\mathbf{u}}\|_{L^{2}}^{2} +\||\nabla \mathbf{d}||\nabla^{2} \mathbf{d}|\|_{L^{2}}^{2}\right)
+\int_{0}^{T} \left( \|\nabla \dot{\mathbf{u}}\|_{L^{2}}^{2} +\||\nabla \mathbf{d}||\Delta\nabla \mathbf{d}| \|_{L^{2}}^{2}\right)\text{d}t\leq C,
\end{equation}
and
\begin{equation}\label{eq3.14}
\sup_{0\leq t\leq T} t^{i}\left( \|\nabla^{2} \mathbf{u}\|_{L^{2}}^{2} +\|\nabla P\|_{L^{2}}^{2}\right)
\leq C.
\end{equation}
\end{lemma}

{\it Proof.}
First, operating $\partial_{t}+\mathbf{u}\cdot\nabla $ to $\eqref{eq1.1}_{2}^j$ ($j=1,2$) yields that
\begin{align*}
&\partial_{t}(\rho \dot{u}_{j})+\divv(\rho \mathbf{u} \dot{u}_{j}) -\Delta \dot{u}_{j}\nonumber\\
& = \Big(-\partial_{i} (\partial_{i} \mathbf{u}\cdot\nabla u_{j})-\divv (\partial_{i} \mathbf{u}\partial_{i}u_{j})  -\partial_{t}\partial_{i}(\partial_{i}\mathbf{d}\partial_{j}\mathbf{d}) -
\mathbf{u}\cdot\nabla \partial_{i}(\partial_{i}\mathbf{d}\partial_{j} \mathbf{d})\Big)-\partial_{t}\partial_{j} P -\mathbf{u}\cdot\nabla\partial_{j} P.
\end{align*}
Now, multiplying the above equality by $\dot{u}_{j}$, and then integrating by parts over $\mathbb{R}^{2}$, we deduce
\begin{align}\label{eq3.15}
\frac12\frac{d}{dt}\int\rho|\dot{\mathbf{u}}|^{2}dx
+\int|\nabla\dot{\mathbf{u}}|^{2}dx
= & - \int \left[\partial_{i}(\partial_{i} \mathbf{u}\cdot\nabla u_{j})+\divv(\partial_{i} \mathbf{u}\partial_{i} u_{j})\right]\dot{u}_{j}\text{d}x\nonumber\\
& -\int (\partial_{t}\partial_{j} P+\mathbf{u}\cdot\nabla\partial_{j} P)\dot{u}_{j}\text{d}x-\int \partial_{t}\partial_{i}
(\partial_{i}\mathbf{d}\partial_{j}\mathbf{d})
\dot{u}_{j}\text{d}x
\nonumber\\
& -\int\mathbf{u}\cdot\nabla\partial_{i}
(\partial_{i}\mathbf{d}\partial_{j}\mathbf{d}) \dot{u}_{j}\text{d}x\nonumber\\
\triangleq & L_{1}+L_{2}+L_{3}+L_{4}.
\end{align}
By the same arguments as \cite[Lemma 3.3]{LSZ}, one has
\begin{align*}
L_{1}+L_{2}\leq \frac{d}{dt} \int P\partial_{j} u_{i}\partial_{i}u_{j}\text{d}x +C(\|P\|_{L^{4}}^{4}+\|\nabla \mathbf{u}\|_{L^{4}}^{4}) +\frac{1}{4} \|\nabla \dot{\mathbf{u}}\|_{L^{2}}^{2}.
\end{align*}
For terms $L_{3}$ and $L_{4}$, we infer from $\eqref{eq1.1}_{3}$ and $\eqref{eq1.1}_{4}$ that
\begin{align*}
L_{3}+L_{4} = & \int \partial_{i}\dot{u}_{j} \partial_{i} \mathbf{d}_{t}\partial_{j}\mathbf{d}\text{d}x +\int \partial_{i}\dot{u}_{j}\partial_{i}\mathbf{d}\partial_{j}\mathbf{d}_{t}
\text{d}x-\int u_{k}\partial_{k}\partial_{i}(\partial_{i}\mathbf{d}\partial_{j}\mathbf{d})
\dot{u}_{j}\text{d}x\nonumber\\
= & \int\! \partial_{i}\dot{u}_{j}  \partial_{j}\mathbf{d}\partial_{i} (\Delta\mathbf{d}-\mathbf{u}\cdot \nabla\mathbf{d}+|\nabla\mathbf{d}|^{2}\mathbf{d})\text{d}x +\int\! \partial_{i}\dot{u}_{j}\partial_{i}\mathbf{d}\partial_{j}(\Delta\mathbf{d}
-\mathbf{u}\cdot \nabla\mathbf{d}+|\nabla\mathbf{d}|^{2}\mathbf{d})\text{d}x\nonumber\\
& -\int u_{k}\partial_{k}\partial_{i}(\partial_{i}\mathbf{d}\partial_{j}\mathbf{d})
\dot{u}_{j}\text{d}x\nonumber\\
= & \int\! \partial_{i}\dot{u}_{j}  \partial_{j}\mathbf{d}(\Delta \partial_{i}\mathbf{d}-\partial_{i}\mathbf{u}\cdot \nabla\mathbf{d}+\partial_{i}(|\nabla \mathbf{d}|^{2}\mathbf{d}))\text{d}x-\int \partial_{i}\dot{u}_{j} \partial_{j}\mathbf{d} u_{k}\partial_{k}\partial_{i}\mathbf{d}\text{d}x\nonumber\\
&+\int\! \partial_{i}\dot{u}_{j}\partial_{i}\mathbf{d}(\Delta\partial_{j} \mathbf{d}-\partial_{j}\mathbf{u}\cdot \nabla\mathbf{d}+\partial_{j}(|\nabla \mathbf{d}|^{2}\mathbf{d}))\text{d}x-\int \partial_{i}\dot{u}_{j} \partial_{i}\mathbf{d} u_{k}\partial_{k}\partial_{j}\mathbf{d}\text{d}x\nonumber\\
&+\int \partial_{i} u_{k}\partial_{k}(\partial_{i}\mathbf{d}\partial_{j}\mathbf{d})
\dot{u}_{j}\text{d}x
+\int u_{k}\partial_{k}
(\partial_{i}\mathbf{d}\partial_{j}\mathbf{d})\partial_{i} \dot{u}_{j}\text{d}x\nonumber\\
=&\int\! \partial_{i}\dot{u}_{j}  \partial_{j}\mathbf{d} (\Delta \partial_{i}\mathbf{d}-\partial_{i}\mathbf{u}\cdot \nabla\mathbf{d}+\partial_{i}(|\nabla\mathbf{d}|^{2}\mathbf{d}))\text{d}x\nonumber\\
&+\!\int\! \partial_{i}\dot{u}_{j}\partial_{i}\mathbf{d}(\Delta\partial_{j} \mathbf{d}-\partial_{j} \mathbf{u}\cdot \nabla \mathbf{d}+\partial_{j}(|\nabla \mathbf{d}|^{2}\mathbf{d}))\text{d}x +\int \partial_{i} u_{k}\partial_{i}\mathbf{d}\partial_{j}\mathbf{d}\partial_{k}\dot{u}_{j}
\text{d}x,
\end{align*}
which combined with $|\mathbf{d}|=1$, H{\"o}lder's and Young's inequalities yields
\begin{align*}
L_{3}+L_{4}\leq & C\|\nabla\dot{\mathbf{u}}\|_{L^{2}} \||\nabla \mathbf{d}||\Delta \nabla \mathbf{d}|\|_{L^{2}} +C\|\nabla \dot{\mathbf{u}}\|_{L^{2}} \|\nabla \mathbf{u}\|_{L^{4}} \||\nabla\mathbf{d}|^{2}\|_{L^{4}}\nonumber\\
& + C\|\nabla \dot{\mathbf{u}}\|_{L^{2}}\||\nabla\mathbf{d}|^{4}\|_{L^{2}}
+C\|\nabla \dot{\mathbf{u}}\|_{L^{2}}\||\nabla\mathbf{d}|^{2}|\nabla^{2}\mathbf{d}|\|_{L^{2}}\nonumber \\
\leq &\frac{1}{4} \|\nabla \dot{\mathbf{u}}\|_{L^{2}}^{2}+C\||\nabla \mathbf{d}||\Delta \nabla \mathbf{d}|\|_{L^{2}}^{2} +C\|\nabla \mathbf{u}\|_{L^{4}}^{4}+C \||\nabla \mathbf{d}|^{2}\|_{L^{4}}^{4}
+C\||\nabla\mathbf{d}|^2|\nabla^2 \mathbf{d}|\|_{L^2}^2.
\end{align*}
Inserting the above estimates of $L_{i}\ (i=1,2,3,4)$ into \eqref{eq3.15}, one deduces
\begin{align}\label{eq3.16}
&\frac{1}{2}\frac{d}{dt} \|\sqrt{\rho}\dot{\mathbf{u}}\|_{L^{2}}^{2}+\frac{1}{2} \|\nabla \dot{\mathbf{u}}\|_{L^{2}}^{2}\nonumber\\
& \leq \frac{d}{dt} \int P\partial_{j} u_{i}\partial_{i}u_{j}\text{d}x
+C(\|P\|_{L^{4}}^{4}+\|\nabla \mathbf{u}\|_{L^{4}}^{4}+\||\nabla \mathbf{d}|^{2}\|_{L^{4}}^{4})\nonumber \\
& \quad +C\||\nabla \mathbf{d}||\Delta \nabla \mathbf{d}|\|_{L^{2}}^{2}+C\||\nabla\mathbf{d}|^2|\nabla^2 \mathbf{d}|\|_{L^2}^2.
\end{align}

Now, inspired by \cite{LH15,LXZ,LXZ1}, for $a_{1},a_{2}\in \{-1,0,1\}$, let us denote
\begin{align*}
(\widetilde{\nabla} \mathbf{d})(a_{1},a_{2})=a_{1}\partial_{1}\mathbf{d}+a_{2}\partial_{2}\mathbf{d}, \quad (\widetilde{\nabla}\mathbf{u})(a_{1},a_{2})=a_{1}\partial_{1} \mathbf{u}+a_{2}\partial_{2}\mathbf{u},\quad \widetilde{u}(a_{1},a_{2})=a_{1}u_{1}+a_{2} u_{2},
\end{align*}
then it is easy to deduce from \eqref{eq3.7} that
\begin{align*}
\widetilde{\nabla }\mathbf{d}_{t} -\Delta\widetilde{\nabla } \mathbf{d}=-\widetilde{\nabla}\mathbf{u}\cdot\nabla\mathbf{d} -\mathbf{u}\cdot\nabla\widetilde{\nabla}\mathbf{d}+|\nabla \mathbf{d}|^{2}\widetilde{\nabla}\mathbf{d}+ 2(\nabla \mathbf{d}\cdot\mathbf{d})\nabla\widetilde{\nabla}\mathbf{d}.
\end{align*}
Multiplying the above equality by $4\widetilde{\nabla}\mathbf{d}\Delta|\widetilde{\nabla} \mathbf{d}|^{2}$, and then integrating by parts over $\mathbb{R}^{2}$, it follows that
\begin{align}\label{eq3.17}
&\frac{d}{dt} \|\nabla |\widetilde{\nabla }\mathbf{d}|^{2}\|_{L^{2}}^{2}+2\|\Delta |\widetilde{\nabla}\mathbf{d}|^{2}\|_{L^{2}}^{2}\nonumber\\
& = -4\int (\widetilde{\nabla}\mathbf{u}\cdot\nabla\mathbf{d})\cdot\widetilde{\nabla} \mathbf{d}\Delta|\widetilde{\nabla}\mathbf{d}|^{2}\text{d}x
-4\int (\mathbf{u}\cdot\nabla \widetilde{\nabla}\mathbf{d})\cdot\widetilde{\nabla}\mathbf{d}
\Delta|\widetilde{\nabla }\mathbf{d}|^{2}\text{d}x\nonumber\\
& \quad +4\int |\nabla\mathbf{d}|^{2}|\widetilde{\nabla}\mathbf{d}|^{2} \Delta|\widetilde{\nabla}\mathbf{d}|^{2}\text{d}x
+8\int(\nabla\mathbf{d}\cdot\mathbf{d})\nabla\widetilde{\nabla} \mathbf{d}\cdot\widetilde{\nabla}\mathbf{d}\Delta|\widetilde{\nabla}\mathbf{d}|^{2}\text{d}x\nonumber\\
& \leq C\int |\widetilde{\nabla}\mathbf{u}| |\nabla \mathbf{d}||\widetilde{\nabla}\mathbf{d}||\Delta|\widetilde{\nabla}\mathbf{d}|^{2}|\text{d}x
-2\int (\mathbf{u}\cdot\nabla |\widetilde{\nabla }\mathbf{d}|^{2})\Delta|\widetilde{\nabla}\mathbf{d}|^{2}\text{d}x\nonumber\\
& \quad +C\int |\nabla\mathbf{d}|^{2}|\widetilde{\nabla}\mathbf{d}|^{2} |\Delta|\widetilde{\nabla}\mathbf{d}|^{2}|\text{d}x+C\int |\nabla\mathbf{d}| \nabla|\widetilde{\nabla} \mathbf{d}|^{2}|\Delta|\widetilde{\nabla}\mathbf{d}|^{2}|\text{d}x\nonumber\\
& \leq C\int |\widetilde{\nabla} \mathbf{u}| |\nabla \mathbf{d}||\widetilde{\nabla}\mathbf{d}||\Delta|\widetilde{\nabla}\mathbf{d}|^{2}|\text{d}x
+C\int |\nabla\mathbf{u}||\nabla |\widetilde{\nabla }\mathbf{d}|^{2}|^{2}\text{d}x\nonumber\\
& \quad +C\int |\nabla\mathbf{d}|^{2}|\widetilde{\nabla}\mathbf{d}|^{2} |\Delta|\widetilde{\nabla}\mathbf{d}|^{2}|\text{d}x+C\int |\nabla \mathbf{d}| \nabla|\widetilde{\nabla} \mathbf{d}|^{2}|\Delta|\widetilde{\nabla}\mathbf{d}|^{2}|\text{d}x\nonumber\\
& \leq C\|\nabla\mathbf{u}\|_{L^{4}}\||\nabla\mathbf{d}|^{2}\|_{L^{4}} \|\Delta|\widetilde{\nabla }\mathbf{d}|^{2}\|_{L^{2}}
+C\|\nabla\mathbf{u}\|_{L^{4}} \|\nabla |\widetilde{\nabla}\mathbf{d}|^{2}\|_{L^{\frac{8}{3}}}^{2}\nonumber\\
& \quad +C\||\nabla\mathbf{d}|^{2}\|_{L^{4}}^{2}
\|\Delta|\widetilde{\nabla}\mathbf{d}|^{2}\|_{L^{2}} +C\|\Delta|\widetilde{\nabla}\mathbf{d}|^{2}\|_{L^{2}} \|\nabla|\widetilde{\nabla}\mathbf{d}|^{2}\|_{L^{\frac{8}{3}}} \|\nabla \mathbf{d}\|_{L^{8}}\nonumber\\
& \leq \frac{1}{2}\|\Delta|\widetilde{\nabla}\mathbf{d}|^{2}\|_{L^{2}}^{2}
+C(\||\nabla \mathbf{d}|^{2}\|_{L^{4}}^{4}+\|\nabla\mathbf{u}\|_{L^{4}}^{4})
+C\|\nabla |\widetilde{\nabla}\mathbf{d}|^{2}\|_{L^{\frac{8}{3}}}^{\frac{8}{3}} +C\|\nabla \mathbf{d}\|_{L^{8}}^{8}\nonumber\\
& \leq \frac{1}{2}\|\Delta|\widetilde{\nabla}\mathbf{d}|^{2}\|_{L^{2}}^{2}+C(\||\nabla \mathbf{d}|^{2}\|_{L^{4}}^{4}+\|\nabla\mathbf{u}\|_{L^{4}}^{4})
+C\||\widetilde{\nabla} \mathbf{d}|^{2}\|_{L^{4}}^{\frac{4}{3}}\|\Delta|\widetilde{\nabla} \mathbf{d}|^{2}\|_{L^{2}}^{\frac{4}{3}} +C\|\nabla\mathbf{d}\|_{L^{2}}^{4}\|\nabla^{2} \mathbf{d}\|_{L^{4}}^{4}\nonumber\\
& \leq \|\Delta|\widetilde{\nabla}\mathbf{d}|^{2}\|_{L^{2}}^{2}+C(\||\nabla \mathbf{d}|^{2}\|_{L^{4}}^{4}+\|\nabla\mathbf{u}\|_{L^{4}}^{4})+C\|\nabla^{2} \mathbf{d}\|_{L^{4}}^{4},
\end{align}
where we have used H\"{o}lder's and  Gagliardo-Nirenberg inequalities, and \eqref{ll} in the above estimates. Noticing that
\begin{align*}
\||\Delta\nabla\mathbf{d}||\nabla\mathbf{d}|\|_{L^{2}}^{2}\!
\leq & \! C\|\nabla^{2}\mathbf{d}\|_{L^{4}}^{4}\!+\! \|\Delta|\widetilde{\nabla }\mathbf{d}(1,0)|^{2}\|_{L^{2}}^{2}\!+\!\|\Delta|\widetilde{\nabla }\mathbf{d}(0,1)|^{2}\|_{L^{2}}^{2}\!+\!\|\Delta|\widetilde{\nabla }\mathbf{d}(1,1)|^{2}\|_{L^{2}}^{2}\! \\
& +\!\|\Delta|\widetilde{\nabla }\mathbf{d}(1,-1)|^{2}\|_{L^{2}}^{2},
\end{align*}
and
\begin{align}\label{eq3.18}
\||\nabla^{2}\mathbf{d}||\nabla\mathbf{d}|\|_{L^{2}}^{2}\leq G(t)\leq C\||\nabla^{2}\mathbf{d}||\nabla\mathbf{d}|\|_{L^{2}}^{2}
\end{align}
with
\begin{align*}
G(t)\triangleq \|\nabla|\widetilde{\nabla}\mathbf{d}(1,0)|^{2}\|_{L^{2}}^{2}
+\|\nabla|\widetilde{\nabla}\mathbf{d}(0,1)|^{2}\|_{L^{2}}^{2}
+\|\nabla|\widetilde{\nabla}\mathbf{d}(1,1)|^{2}\|_{L^{2}}^{2}
+\|\nabla|\widetilde{\nabla}\mathbf{d}(1,-1)|^{2}\|_{L^{2}}^{2}.
\end{align*}
Thus, it follows from \eqref{eq3.17} multiplied by $(C_{2}+1)$ that
\begin{align*}
\frac{d}{dt} \left((C_{2}+1)G(t)\right) +(C_{2}+1)\||\Delta\nabla\mathbf{d}||\nabla \mathbf{d}|\|_{L^{2}}^{2}
\leq  C\|\nabla\mathbf{u}\|_{L^{4}}^{4} +C\|\nabla^{2}\mathbf{d}\|_{L^{4}}^{4}+C\||\nabla\mathbf{d}|^{2}\|_{L^{4}}^{4},
\end{align*}
which combined with \eqref{eq3.16} ensures that
\begin{align}\label{eq3.19}
\frac{d}{dt}F(t)+\frac{1}{2}\|\nabla \dot{\mathbf{u}}\|_{L^{2}}^{2}+\||\Delta\nabla \mathbf{d}||\nabla \mathbf{d}|\|_{L^{2}}^{2}
\leq C\|P\|_{L^{4}}^{4}+C\|\nabla\mathbf{u}\|_{L^{4}} +C\|\nabla^{2} \mathbf{d}\|_{L^{4}}^{4}+C\||\nabla \mathbf{d}|^{2}\|_{L^{4}}^{4},
\end{align}
where
\begin{align*}
F(t)\triangleq \frac{1}{2}\|\sqrt{\rho}\dot{\mathbf{u}}\|_{L^{2}}^{2}+(C_{2}+1) G(t)-\int P\partial_{j}u_{i}\partial_{i}u_{j}\text{d}x
\end{align*}
satisfies
\begin{align}\label{eq3.20}
\frac{1}{4}\|\sqrt{\rho}\dot{\mathbf{u}}\|_{L^{2}}^{2}+\frac{C_{2}+1}{2} G(t)-C\|\nabla\mathbf{u}\|_{L^{4}}^{4}
\leq F(t)
\leq \|\sqrt{\rho}\dot{\mathbf{u}}\|_{L^{2}}^{2}+C G(t)+C\|\nabla \mathbf{u}\|_{L^{4}}^{4}
\end{align}
owing to the following estimate
\begin{align*}
\left|\int P\partial_{j}u_{i}\partial_{i}u_{j}\text{d}x\right|
\leq & C\|P\|_{BMO} \|\partial_{i}u_{j}\partial_{j} u_{i}\|_{\mathcal{H}^{1}}\leq C\|\nabla P\|_{L^{2}}\|\nabla \mathbf{u}\|_{L^{2}}^{2}\quad (\text{by Lemma \ref{lem2.5}})\nonumber\\
\leq &C(\|\sqrt{\rho}\dot{\mathbf{u}}\|_{L^{2}} +\||\nabla^{2}\mathbf{d}||\nabla \mathbf{d}|\|_{L^{2}})\|\nabla\mathbf{u}\|_{L^{2}}^{2}
\quad (\text{by} \eqref{eq3.9})\nonumber\\
\leq& \frac{1}{2}\|\sqrt{\rho}\dot{\mathbf{u}}\|_{L^{2}}^{2} +\frac{C_{2}+1}{2} G(t)+C\|\nabla\mathbf{u}\|_{L^{4}}^{4}.
\end{align*}

Next, we shall estimate the terms on the right-hand side of \eqref{eq3.19}. To bound the terms $\|P\|_{L^{4}}$ and $\|\nabla u\|_{L^{4}}$, it follows from Sobolev embedding, \eqref{eq3.9}, H\"{o}lder's inequality, \eqref{ll}, \eqref{eq3.5}, and \eqref{eq3.20}  that
\begin{align}\label{eq3.21}
\|P\|_{L^{4}}^{4} +\|\nabla\mathbf{u}\|_{L^{4}}^{4}
\leq & C(\|\nabla P\|_{L^{\frac{4}{3}}}^{4}
+\|\nabla^{2}\mathbf{u}\|_{L^{\frac{4}{3}}}^{4})
\leq C(\|\rho\dot{\mathbf{u}}\|_{L^{\frac{4}{3}}}^{4}
+\||\nabla\mathbf{d}||\nabla^{2} \mathbf{d}|\|_{L^{\frac{4}{3}}}^{4})\nonumber\\
\leq & C\|\rho\|_{L^{2}}^{2}\|\sqrt{\rho} \dot{\mathbf{u}}\|_{L^{2}}^{4}+C\|\nabla\mathbf{d}\|_{L^{2}}^{4}\|\nabla^{2} \mathbf{d}\|_{L^{4}}^{4}\nonumber\\
\leq& C \|\sqrt{\rho} \dot{\mathbf{u}}\|_{L^{2}}^{2}(F(t)+\|\nabla \mathbf{u}\|_{L^{2}}^{4})
+C\|\nabla^{2}\mathbf{d}\|_{L^{2}}^{2}\|\nabla^{3}\mathbf{d}\|_{L^{2}}^{2}.
\end{align}
By Ladyzhenskaya's inequality, \eqref{ll}, and \eqref{eq3.20}, one has
\begin{align}\label{eq3.22}
\|\nabla^{2}\mathbf{d}\|_{L^{4}}^{4}+\||\nabla\mathbf{d}|^{2}\|_{L^{4}}^{4}
\leq & C\|\nabla^{2}\mathbf{d}\|_{L^{2}}^{2}\|\nabla^{3}\mathbf{d}\|_{L^{2}}^{2}
+C\|\nabla \mathbf{d}\|_{L^{2}}^{2} \|\nabla^{2}\mathbf{d}\|_{L^{2}}^{2} \||\nabla \mathbf{d}||\nabla^{2}\mathbf{d}|\|_{L^{2}}^{2}\nonumber\\
\leq & C\|\nabla^{2}\mathbf{d}\|_{L^{2}}^{2}(F(t)+\|\nabla \mathbf{u}\|_{L^{4}}^{4})
+C\|\nabla^{2}\mathbf{d}\|_{L^{2}}^{2}\|\nabla^{3}\mathbf{d}\|_{L^{2}}^{2}.
\end{align}
Then, substituting \eqref{eq3.21} and \eqref{eq3.22} into \eqref{eq3.19}, one obtains
\begin{align}\label{eq3.23}
&\frac{d}{dt} F(t) +\frac{1}{2}\|\nabla \dot{\mathbf{u}}\|_{L^{2}}^{2}+\||\Delta\nabla \mathbf{d}||\nabla \mathbf{d}|\|_{L^{2}}^{2}\nonumber\\
& \leq C(\|\sqrt{\rho} \dot{\mathbf{u}}\|_{L^{2}}^{2}+\|\nabla^{2} \mathbf{d}\|_{L^{2}}^{2})(F(t)+\|\nabla\mathbf{u}\|_{L^{2}}^{4})
+C\|\nabla^{2}\mathbf{d}\|_{L^{2}}^{2}\|\nabla^{3}\mathbf{d}\|_{L^{2}}^{2}.
\end{align}
Multiplying \eqref{eq3.23} by $t^{i}\ (i=1,2)$, and then applying Gronwall's inequality, it follows from \eqref{eq3.18}, \eqref{eq3.20}, \eqref{eq3.2}, \eqref{eq3.3}, and \eqref{ll} that
\begin{align*}
&\sup_{0\leq t\leq T} (t^{i}F(t))+\int_{0}^{T} t^{i}(\|\nabla \dot{\mathbf{u}}\|_{L^{2}}^{2}+\||\Delta\nabla\mathbf{d}||\nabla \mathbf{d}|\|_{L^{2}}^{2})\text{d}t\nonumber\\
& \leq C\int_{0}^{T} t^{i-1} F(t)\text{d}t+C\int_{0}^{T} t^{i} \|\nabla^{2} \mathbf{d}\|_{L^{2}}^{2}\|\nabla^{3}\mathbf{d}\|_{L^{2}}^{2}\text{d}t
+C\int_{0}^{T} (\|\sqrt{\rho} \dot{\mathbf{u}}\|_{L^{2}}^{2}+\|\nabla^{2} \mathbf{d}\|_{L^{2}}^{2})t^{i} \|\nabla\mathbf{u}\|_{L^{4}}^{4}\text{d}t\nonumber\\
& \leq C\int_{0}^{T} t^{i-1} (\|\sqrt{\rho} \dot{\mathbf{u}}\|_{L^{2}}^{2}+\||\nabla^{2}\mathbf{d}||\nabla \mathbf{d}|\|_{L^{2}}^{2})\text{d}t+C\sup_{0\leq t\leq T}(t^{i-1}\|\nabla \mathbf{u}\|_{L^{2}}^{2})\int_{0}^{T} \|\nabla \mathbf{u}\|_{L^{2}}^{2}\text{d}t\nonumber\\
& \quad +C\!\sup_{0\leq t\leq T} (t^{i-1}\|\nabla^{2}\mathbf{d}\|_{L^{2}}^{2})\int_{0}^{T} t\|\nabla^{3}\mathbf{d}\|_{L^{2}}^{2}\text{d}t
+C\!\sup_{0\leq t\leq T}(t^{i} \|\nabla\mathbf{u}\|_{L^{4}}^{4})\int_{0}^{T} (\|\sqrt{\rho} \dot{\mathbf{u}}\|_{L^{2}}^{2}+\|\nabla^{2} \mathbf{d}\|_{L^{2}}^{2})\text{d}t\nonumber\\
& \leq C\int_{0}^{T} t^{i-1} (\|\sqrt{\rho} \dot{\mathbf{u}}\|_{L^{2}}^{2}+\|\nabla^{2}\mathbf{d}\|_{L^{4}}^{2}\|\nabla \mathbf{d}\|_{L^{4}}^{2})\text{d}t+C\nonumber\\
& \leq C\int_{0}^{T} t^{i-1} (\|\sqrt{\rho} \dot{\mathbf{u}}\|_{L^{2}}^{2}+ \|\nabla\mathbf{d}\|_{L^{2}}\|\nabla^{2}\mathbf{d}\|_{L^{2}}^{2} \|\nabla^{3}\mathbf{d}\|_{L^{2}})\text{d}t+C\nonumber\\
& \leq C\int_{0}^{T} t^{i-1} (\|\sqrt{\rho} \dot{\mathbf{u}}\|_{L^{2}}^{2}+ \|\nabla^{3}\mathbf{d}\|_{L^{2}}^{2})\text{d}t+C\sup_{0\leq t\leq T} (t^{i-1}\|\nabla^{2}\mathbf{d}\|_{L^{2}}^{2})\int_{0}^{T} \|\nabla^{2}\mathbf{d}\|_{L^{2}}^{2}\text{d}t+C\nonumber\\
& \leq C.
\end{align*}
This together with \eqref{eq3.18}, \eqref{eq3.20}, \eqref{eq3.2}, \eqref{eq3.3}, and \eqref{ll} implies \eqref{eq3.13}.

 Finally, it is easy to see that the estimate \eqref{eq3.13} combined with \eqref{eq3.9} gives \eqref{eq3.14}. This completes the proof of Lemma \ref{lem3.3}.   \hfill $\Box$

\subsection{Higher order estimates}
The following spatial weighted estimate on the density plays an important role in deriving the bounds on the higher order derivatives of the solutions $(\rho,\mathbf{u},P,\mathbf{d})$, whose proof can be found in \cite[Lemma 3.4]{LSZ}.
\begin{lemma}\label{lem3.4}
There exists a positive constant $C$ depending on $T$ such that
\begin{align}\label{eq3.24}
\sup_{0\leq t\leq T} \|\rho\bar{x}^{a}\|_{L^{1}}\leq C(T).
\end{align}
\end{lemma}

\begin{lemma}\label{lem3.5}
There exists a positive constant $C$ depending on $T$ such that
\begin{align}\label{eq3.26}
\sup_{t\in[0,T]}\|\rho\|_{H^{1}\cap W^{1,q}}
&+\int_{0}^{T}\left(\|\nabla^{2} \mathbf{u}\|_{L^2}^2+\|\nabla^{2} \mathbf{u} \|_{L^q}^{\frac{q+1}{q}}
+t\|\nabla^{2} \mathbf{u}\|_{L^2 \cap L^q}^2\right)dt \notag \\
&+\int_{0}^{T}\left(\|\nabla P\|_{L^2}^2+\|\nabla P \|_{L^q}^\frac{q+1}{q}+t\|\nabla P\|_{L^2\cap L^q}^2\right)dt
\leq C(T).
\end{align}
\end{lemma}

{\it Proof.}
First, it follows from $\eqref{eq1.1}_{1}$ and $\eqref{eq1.1}_{4}$ that
$\nabla\rho$ satisfies for any $r\geq2$,
\begin{align}\label{eq3.27}
\frac{d}{dt} \|\nabla \rho\|_{L^{r}}\leq C\|\nabla \mathbf{u}\|_{L^{\infty}} \|\nabla \rho\|_{L^{r}}.
\end{align}
Next, by Gagliardo-Nirenberg inequality, \eqref{eq3.2}, and \eqref{eq3.9}, one gets for $q>2$ as in Theorem \ref{thm1.2},
\begin{align}\label{eq3.28}
\|\nabla \mathbf{u}\|_{L^{\infty}}
\leq C\|\nabla \mathbf{u}\|_{L^{2}}^{\frac{q-2}{2(q-1)}} \|\nabla^{2} \mathbf{u}\|_{L^{q}}^{\frac{q}{2(q-1)}}
\leq C\left( \|\rho\dot{\mathbf{u}}\|_{L^{q}}^{\frac{q}{2(q-1)}}+\||\nabla \mathbf{d}||\nabla^{2}\mathbf{d}|\|_{L^{q}}^{\frac{q}{2(q-1)}}\right).
\end{align}
By virtue of Gagliardo-Nirenberg inequality and Lemma \ref{lem2.4}, one has
\begin{align*}
\|\rho\dot{\mathbf{u}}\|_{L^{q}}
\leq  & C\|\rho \dot{\mathbf{u}}\|_{L^{2}}^{\frac{2q^2-1}{q(q^2-1)}} \|\rho\dot{\mathbf{u}}\|_{L^{2q^2}}^{\frac{q^2-2q}{q^2-1}}
\nonumber\\
\leq &  C\|\sqrt{\rho}
\dot{\mathbf{u}}\|_{L^{2}}^{\frac{2q^2-1}{q(q^2-1)}}
(\|\sqrt{\rho} \dot{\mathbf{u}}\|_{L^{2}}
+\|\nabla \dot{\mathbf{u}}\|_{L^{2}})^{\frac{q^2-2q}{q^2-1}}
\nonumber\\
\leq & C \|\sqrt{\rho} \dot{\mathbf{u}}\|_{L^{2}}
+C\|\sqrt{\rho}\dot{\mathbf{u}}\|_{L^{2}}^{\frac{2q^2-1}{q(q^2-1)}}
\|\nabla \dot{\mathbf{u}}\|_{L^{2}}^{\frac{q^2-2q}{q^2-1}},
\end{align*}
which along with \eqref{eq3.2} and \eqref{eq3.13} leads to
\begin{align}\label{eq3.29}
&\int_{0}^{T}\left(\|\rho\dot{\mathbf{u}}\|_{L^{q}}^{\frac{q+1}{q}}
+t\|\rho\dot{\mathbf{u}}\|_{L^{q}}^{2}\right)\text{d}t\nonumber\\
& \leq C\int_{0}^{T} t^{-\frac{q+1}{2q}} \left(t\|\sqrt{\rho}\dot{\mathbf{u}}\|_{L^{2}}^{2}\right)^{\frac{2q-1}{2q(q-1)}}
\left(t\|\nabla \dot{\mathbf{u}}\|_{L^{2}}^{2}\right)^{\frac{q-2}{2q-2}}\text{d}t +C\int_{0}^{T} \|\sqrt{\rho} \dot{\mathbf{u}}\|_{L^{2}}^{\frac{q+1}{q}}\text{d}t\nonumber\\
& \quad +C\int_{0}^{T} \left(t\|\sqrt{\rho} \dot{\mathbf{u}}\|_{L^{2}}^{2}\right)^{\frac{2q-1}{q^2-1}} \left( t \|\nabla \dot{\mathbf{u}}\|_{L^{2}}^{2}\right)^{\frac{q(q-2)}{q^2-1}}\text{d}t
+C\int_{0}^{T} t\|\sqrt{\rho} \dot{\mathbf{u}}\|_{L^{2}}^{2}\text{d}t\nonumber\\
& \leq C\sup_{0\leq t\leq T}\left(t\|\sqrt{\rho}\dot{\mathbf{u}}\|_{L^{2}}^{2}\right)
^{\frac{2q-1}{2q(q-1)}} \int_{0}^{T}\left(t\|\nabla \dot{\mathbf{u}}\|_{L^{2}}^{2}+t^{-\frac{q^{3}+q^{2}-q-1}{q^{3}+q^{2}}}\right)
\text{d}t\nonumber\\
& \quad +C\int_{0}^{T}\left(t\|\sqrt{\rho} \dot{\mathbf{u}}\|_{L^{2}}^{2}+t\|\nabla \dot{\mathbf{u}}\|_{L^{2}}^{2}\right)\text{d}t+C\int_{0}^{T}
(1+\|\sqrt{\rho} \dot{\mathbf{u}}\|_{L^{2}}^{2})\text{d}t\nonumber\\
& \leq  C(T).
\end{align}
On the other hand, it follows from H\"{o}lder's and Gagliardo-Nirenberg inequalities,  and \eqref{eq3.2} that
\begin{align}\label{eq3.30}
&\int_{0}^{T} \left(\||\nabla \mathbf{d}||\nabla^{2}\mathbf{d}|\|_{L^{q}}^{\frac{q+1}{q}}+t\||\nabla \mathbf{d}||\nabla^{2}\mathbf{d}|\|_{L^{q}}^{2}\right)\text{d}t\nonumber\\
& \leq C\int_{0}^{T} \left[\left(\|\nabla \mathbf{d}\|_{L^{2q}}\|\nabla^{2}\mathbf{d}\|_{L^{2q}}\right)^{\frac{q+1}{q}}
+t\left(\|\nabla \mathbf{d}\|_{L^{2q}}\|\nabla^{2}\mathbf{d}\|_{L^{2q}}\right)^{2}\right]
\text{d}t\nonumber\\
& \leq C\int_{0}^{T}  \left[\left(\|\nabla\mathbf{d}\|_{L^{2}}^{\frac{1}{q}}\|\nabla^{2} \mathbf{d}\|_{L^{2}}\|\nabla^{3}\mathbf{d}\|_{L^{2}}^{1-\frac{1}{q}}\right)
^{\frac{q+1}{q}}+t\left(\|\nabla \mathbf{d}\|_{L^{2}}^{\frac{1}{q}}\|\nabla^{2} \mathbf{d}\|_{L^{2}}\|\nabla^{3}\mathbf{d}\|_{L^{2}}
^{1-\frac{1}{q}}\right)^{2}\right]\text{d}t\nonumber\\
& \leq C\int_{0}^{T} \left(\|\nabla^{3} \mathbf{d}\|_{L^{2}}^{2}+t^{q}+1\right)\text{d}t\nonumber\\
& \leq C(T).
\end{align}
Hence, combining \eqref{eq3.28}, \eqref{eq3.29}, and \eqref{eq3.30} together, it follows that
\begin{align}\label{eq3.31}
\int_{0}^{T} \|\nabla\mathbf{u}\|_{L^{\infty}}\text{d}t\leq C(T).
\end{align}
Thus, applying Gronwall's inequality to \eqref{eq3.27} ensures
\begin{align}\label{eq3.32}
\sup_{0\leq t\leq T} \|\nabla \rho\|_{L^{2}\cap L^{q}}\leq C(T).
\end{align}

Finally, it is easy to deduce from \eqref{eq3.9}, \eqref{eq3.29}, \eqref{eq3.30}, \eqref{eq3.2}, and \eqref{eq3.4} that
\begin{align*}
\int_{0}^{T} \!\!\left(\|\nabla^{2}\mathbf{u}\|_{L^{2}}^{2}\!+\!\|\nabla P\|_{L^{2}}^{2}\!+\!\|\nabla^{2}\mathbf{u}\|_{L^{q}}^{\frac{q+1}{q}}\!\!
+\!\|\nabla P\|_{L^{q}}^{\frac{q+1}{q}}\!\!+t(\|\nabla^{2}\mathbf{u}\|_{L^{2}\cap L^{q}}^2\!+\!\|\nabla P\|_{L^{2}\cap L^{q}}^2)\!\right)\text{d}t
\leq C(T).
\end{align*}
This together with \eqref{eq3.5} and \eqref{eq3.32} yields \eqref{eq3.26}, and finishes the proof of Lemma \ref{lem3.5}.  \hfill $\Box$

We shall now give some spatial estimates on $\nabla \rho$, $\nabla\mathbf{d}$ and $\nabla^{2}\mathbf{d}$, which are crucial to derive the estimates on the gradients of both $\mathbf{u}_{t}$ and $\nabla\mathbf{d}_{t}$.
\begin{lemma}\label{lem3.6}
There exists a positive constant $C$ depending on $T$ such that
\begin{align}\label{eq3.33}
&\sup_{0\leq t\leq T} \|\rho \bar{x}^{a}\|_{L^{1}\cap H^{1}\cap W^{1,q}}\leq C(T),\\
      \label{eq3.34}
&\sup_{0\leq t\leq T} \|\nabla\mathbf{d}\bar{x}^{\frac{a}{2}}\|_{L^{2}}^{2} +\int_{0}^{T} \|\nabla^{2} \mathbf{d}\bar{x}^{\frac{a}{2}}\|_{L^{2}}^{2}\text{d}t
\leq C(T),
\end{align}
and
\begin{align}\label{eq3.35}
\sup_{0\leq t\leq T} t\|\nabla^{2}\mathbf{d}\bar{x}^{\frac{a}{2}}\|_{L^{2}}^{2} +\int_{0}^{T} t\|\nabla^{3}\mathbf{d}\bar{x}^{\frac{a}{2}}\|_{L^{2}}^{2}\text{d}t
\leq C(T).
\end{align}
\end{lemma}

{\it Proof.}
With \eqref{eq3.24} in hand, the proof of \eqref{eq3.33} is exactly the same as \cite[Lemma 3.6]{LSZ}, and we omit it for simplicity. To prove
\eqref{eq3.34}, by multiplying \eqref{eq3.7} with $\nabla\mathbf{d}\bar{x}^{a}$ and integrating by parts yield
\begin{align}\label{eq3.36}
&\frac{1}{2}\frac{d}{dt} \|\nabla \mathbf{d}\bar{x}^{\frac{a}{2}}\|_{L^{2}}^{2} +\|\nabla^{2}\mathbf{d}\bar{x}^{\frac{a}{2}}\|_{L^{2}}^{2}\nonumber\\
& \leq \! C\!\!\int \!\!|\nabla\mathbf{d}| |\nabla^{2}\mathbf{d}| \nabla \bar{x}^{a} \text{d}x\! +\!\!\int \! \! |\nabla\mathbf{u}||\nabla\mathbf{d}|^{2} \bar{x}^{a}
\text{d}x +\!\!\int\! |\mathbf{u}||\nabla\mathbf{d}|^{2} \nabla \bar{x}^{a}\text{d}x +\!\!\int\! \! |\nabla\mathbf{d}|^{2} |\nabla^{2}\mathbf{d}| \bar{x}^{a}\text{d}x
\!+\!\!\int\! |\nabla\mathbf{d}|^{3} \nabla\bar{x}^{a}\text{d}x\nonumber\\
& \triangleq J_{1}+J_{2}+J_{3}+J_{4}+J_{5}.
\end{align}
By virtue of H\"{o}lder's and Ladyzhenskaya's inequalities, \eqref{eq3.2}, \eqref{ll}, and \eqref{l3},  we have
\begin{align*}
J_{1}\leq & C\int |\nabla\mathbf{d}| |\nabla^{2}\mathbf{d}| \bar{x}^{a}\text{d}x\leq \frac{1}{10} \|\nabla^{2}\mathbf{d}\bar{x}^{\frac{a}{2}}\|_{L^{2}}^{2}
+C\|\nabla\mathbf{d}\bar{x}^{\frac{a}{2}}\|_{L^{2}}^{2};\\
J_{2}\leq & \|\nabla\mathbf{u}\|_{L^{2}} \|\nabla\mathbf{d} \bar{x}^{\frac{a}{2}}\|_{L^{4}}^{2}\\
\leq & C\|\nabla \mathbf{d}\bar{x}^{\frac{a}{2}}\|_{L^{2}} \left(\|\nabla^{2}\mathbf{d} \bar{x}^{\frac{a}{2}}\|_{L^{2}}+\|\nabla\mathbf{d}\nabla \bar{x}^{\frac{a}{2}}\|_{L^{2}}\right)\\
\leq & \frac{1}{10} \|\nabla^{2}\mathbf{d}\bar{x}^{\frac{a}{2}}\|_{L^{2}}^{2}+C
 \|\nabla\mathbf{d}\bar{x}^{\frac{a}{2}}\|_{L^{2}}^{2};\\
J_{3}\leq & C\int |\mathbf{u}||\nabla\mathbf{d}|^{2}\bar{x}^{a-\frac{3}{4}}\text{d}x
\leq C\|\nabla \mathbf{d}\bar{x}^{\frac{a}{2}}\|_{L^{4}} \|\nabla\mathbf{d} \bar{x}^{\frac{a}{2}}\|_{L^{2}}\|\mathbf{u}\bar{x}^{-\frac{3}{4}}\|_{L^{4}}\\
\leq  &  C\|\nabla \mathbf{d}\bar{x}^{\frac{a}{2}}\|_{L^{4}}^{2} +
C(\|\sqrt{\rho}\mathbf{u}\|_{L^{2}}^{2}
+\|\nabla\mathbf{u}\|_{L^{2}}^{2})\|\nabla\mathbf{d} \bar{x}^{\frac{a}{2}}\|_{L^{2}}^{2}\nonumber\\
\leq  &  \frac{1}{10} \|\nabla^{2}\mathbf{d} \bar{x}^{\frac{a}{2}}\|_{L^{2}}^{2}+C\|\nabla \mathbf{d}\bar{x}^{\frac{a}{2}}\|_{L^{2}}^{2};\\
J_{4}\leq & \|\nabla^{2}\mathbf{d}\bar{x}^{\frac{a}{2}}\|_{L^{2}}
\|\nabla\mathbf{d}\bar{x}^{\frac{a}{2}}\|_{L^{4}} \|\nabla\mathbf{d}\|_{L^{4}}\\
\leq & \frac{1}{20}\|\nabla^{2}\mathbf{d} \bar{x}^{\frac{a}{2}}\|_{L^{2}}^{2}
+ C\|\nabla\mathbf{d}\|_{L^{4}}^{2}\|\nabla \mathbf{d}\bar{x}^{\frac{a}{2}}\|_{L^{4}}^{2}\\
\leq & \frac{1}{20}\|\nabla^{2}\mathbf{d}\bar{x}^{\frac{a}{2}}\|_{L^{2}}^{2}
+ C\left(\|\nabla\mathbf{d}\|_{L^{2}}^{2}
+\|\nabla^{2}\mathbf{d}\|_{L^{2}}^{2}\right)\|\nabla \mathbf{d}\bar{x}^{\frac{a}{2}}\|_{L^{2}}\left(\|\nabla^{2}\mathbf{d} \bar{x}^{\frac{a}{2}}\|_{L^{2}}+\|\nabla\mathbf{d}\nabla \bar{x}^{\frac{a}{2}}\|_{L^{2}}\right)\\
\leq & \frac{1}{10}\|\nabla^{2}\mathbf{d}\bar{x}^{\frac{a}{2}}\|_{L^{2}}^{2}
+ C \|\nabla\mathbf{d}\bar{x}^{\frac{a}{2}}\|_{L^{2}}^{2};\\
J_{5}\leq&\int |\nabla\mathbf{d}|^{3}\bar{x}^{a}\text{d}x\leq \|\nabla \mathbf{d}\|_{L^{2}}\|\nabla\mathbf{d}\bar{x}^{\frac{a}{2}}\|_{L^{4}}^{2}
\leq \frac{1}{10}\|\nabla^{2}\mathbf{d}\bar{x}^{\frac{a}{2}}\|_{L^{2}}^{2}
+ C \|\nabla\mathbf{d}\bar{x}^{\frac{a}{2}}\|_{L^{2}}^{2}.
\end{align*}
Then, inserting the estimates of $J_{i}\ (i=1,2,\cdots,5)$ into \eqref{eq3.36}, we obtain \eqref{eq3.34} after by using Gronwall's inequality.

It remains to show \eqref{eq3.35}.
Multiplying \eqref{eq3.7} by $\Delta \nabla \mathbf{d}\bar{x}^{a}$ and integrating by parts lead to
 \begin{align}\label{eq3.37}
 &\frac{1}{2}\frac{d}{dt} \|\nabla^{2}\mathbf{d}\bar{x}^{\frac{a}{2}}\|_{L^{2}}^{2}+\! \|\nabla^{3} \mathbf{d}\bar{x}^{\frac{a}{2}}\|_{L^{2}}^{2}\nonumber\\
 & =-\!\int\! \nabla \mathbf{d}_{t}\nabla^{2}\mathbf{d}\nabla \bar{x}^{a}\text{d}x+\!\int \!\nabla(\mathbf{u}\cdot\nabla\mathbf{d})\nabla \Delta \mathbf{d}\bar{x}^{a}\text{d}x
-\!\int\! \nabla(|\nabla\mathbf{d}|^{2}\mathbf{d})\nabla \Delta \mathbf{d}\bar{x}^{a}\text{d}x\nonumber\\
 & \quad -2\int \nabla^{3}\mathbf{d}\nabla^{2}\mathbf{d}\nabla\bar{x}^{a}\text{d}x -\int |\nabla^{2} \mathbf{d}|^{2} \nabla^{2}\bar{x}^{a}\text{d}x
 \nonumber\\
 & \leq \frac{1}{4} \!\|\nabla^{3} \mathbf{d}\bar{x}^{\frac{a}{2}}\|_{L^{2}}^{2}\!+\!\int\!|\Delta\nabla \mathbf{d}||\nabla^{2} \mathbf{d}|\nabla\bar{x}^{a}\text{d}x
 \!+\!\int\! |\mathbf{u}||\nabla\mathbf{d}| |\nabla^{3} \mathbf{d}|\nabla\bar{x}^{a}\text{d}x\!+\!\int\! |\mathbf{u}||\nabla\mathbf{d}| |\nabla^{2} \mathbf{d}|\nabla^{2}\bar{x}^{a}\text{d}x\nonumber\\
 & \quad +\!\int\! |\mathbf{u}||\nabla^{2}\mathbf{d}|^{2}\nabla\bar{x}^{a}\text{d}x+\int |\nabla \mathbf{d}|^{3}|\nabla^{2}\mathbf{d}| \nabla \bar{x}^{a}\text{d}x+\int |\nabla \mathbf{d}||\nabla^{2} \mathbf{d}|^{2}\nabla \bar{x}^{a}\text{d}x\nonumber\\
 &\quad +\int |\nabla\mathbf{u}|^{2}|\nabla\mathbf{d}|^{2}\bar{x}^{a}\text{d}x+\int |\nabla \mathbf{d}|^{6}\bar{x}^{a}\text{d}x+\int
 |\nabla \mathbf{d}|^{2}|\nabla^{2}\mathbf{d}|^{2}\bar{x}^{a}\text{d}x+\int |\nabla^{2}\mathbf{d}|^{2} \nabla^{2}\bar{x}^{a}\text{d}x\nonumber\\
 & \triangleq\frac{1}{4} \|\nabla^{3} \mathbf{d}\bar{x}^{\frac{a}{2}}\|_{L^{2}}^{2}+K_{1}+K_{2}+\cdots+K_{10}.
 \end{align}
Using H\"{o}lder's and Ladyzhenskaya's inequalities, \eqref{l3}, \eqref{eq3.2}, \eqref{ll}, and \eqref{eq3.34}, we get
 \begin{align*}
 K_{1}\leq & C\int|\nabla\Delta\mathbf{d}| |\nabla^{2}\mathbf{d}| \bar{x}^{a}\text{d}x\leq \frac{1}{20} \|\nabla^{3}\mathbf{d}\bar{x}^{\frac{a}{2}}\|_{L^{2}}^{2}
 +C\|\nabla^{2} \mathbf{d}\bar{x}^{\frac{a}{2}}\|_{L^{2}}^{2};\nonumber\\
 K_{2}\leq & C\int |\mathbf{u}||\nabla\mathbf{d}||\nabla^{3} \mathbf{d}|\bar{x}^{a-\frac{3}{4}}\text{d}x\leq \frac{1}{40}\|\nabla^{3}\mathbf{d}\bar{x}^{\frac{a}{2}}\|_{L^{2}}^{2}
 +C\|\mathbf{u}\bar{x}^{-\frac{3}{4}}
 \|_{L^{2}}^{2}\|\nabla^{2} \mathbf{d}\bar{x}^{\frac{a}{2}}\|_{L^{4}}^{2}\nonumber\\
 \leq & \frac{1}{40}\|\nabla^{3}\mathbf{d}\bar{x}^{\frac{a}{2}}\|_{L^{2}}^{2}\!
 +\!C\left(\|\sqrt{\rho} \mathbf{u}\|_{\!L^{2}}^{2}\!+\!\|\nabla\mathbf{u}\|_{\!L^{2}}^{2}\right)\|\nabla^{2} \mathbf{d}\bar{x}^{\frac{a}{2}}\|_{\!L^{2}}\left(\|\nabla^{3}\mathbf{d}\bar{x}^{\frac{a}{2}}
 \|_{\!L^{2}}\!+\!\|\nabla^{2}\mathbf{d}\nabla \bar{x}^{\frac{a}{2}}\|_{\!L^{2}}\!\right)\nonumber\\
 \leq & \frac{1}{20}\|\nabla^{3}\mathbf{d}\bar{x}^{\frac{a}{2}}\|_{L^{2}}^{2}
 +C\|\nabla^{2}\mathbf{d}\bar{x}^{\frac{a}{2}}\|_{L^{2}}^{2};\nonumber\\
 K_{3} \leq  &  C\|\nabla^{2}\mathbf{d}\bar{x}^{\frac{a}{2}}\|_{L^{2}}
 \|\mathbf{u}\bar{x}^{-\frac{3}{4}}\|_{L^{4}}
 \|\nabla\mathbf{d}\bar{x}^{\frac{a}{2}}\|_{L^{4}} \nonumber\\
 \leq  &  C\|\nabla^{2}\mathbf{d}\bar{x}^{\frac{a}{2}}\|_{L^{2}}^2+C
 \left(\|\sqrt{\rho}\mathbf{u}\|_{L^{2}}
 +\|\nabla\mathbf{u}\|_{L^{2}}\right)
 \|\nabla\mathbf{d}\bar{x}^{\frac{a}{2}}\|_{L^{2}}\left(\|\nabla^{2} \mathbf{d}\bar{x}^{\frac{a}{2}}\|_{L^{2}}+\|\nabla\mathbf{d} \nabla\bar{x}^{\frac{a}{2}}\|_{L^{2}}\right)\nonumber\\
 \leq  &  C(1+\|\nabla^{2}\mathbf{d}\bar{x}^{\frac{a}{2}}\|_{L^{2}}^{2});\\
 K_{4}\leq & C \|\nabla^{2}\mathbf{d}\bar{x}^{\frac{a}{2}}\|_{L^{2}} \|\mathbf{u}\bar{x}^{-\frac{3}{4}}\|_{L^{4}}
 \|\nabla^{2}\mathbf{d}\bar{x}^{\frac{a}{2}}\|_{L^{4}}
 \nonumber\\
 \leq & C\|\nabla^{2}\mathbf{d}\bar{x}^{\frac{a}{2}}\|_{L^{2}}^{2}
 +C\left(\|\sqrt{\rho} \mathbf{u}\|_{L^{2}}
 +\|\nabla\mathbf{u}\|_{L^{2}}\right)
 \|\nabla^{2}\mathbf{d}\bar{x}^{\frac{a}{2}}\|_{L^{2}}
 \left(\|\nabla^{3} \mathbf{d}\bar{x}^{\frac{a}{2}}\|_{L^{2}}+\|\nabla^{2}\mathbf{d}\nabla \bar{x}^{\frac{a}{2}}\|_{L^{2}}\right)\nonumber\\
 \leq&\frac{1}{20}\|\nabla^{3}\mathbf{d}\bar{x}^{\frac{a}{2}}\|_{L^{2}}^{2}
 +C\|\nabla^{2}\mathbf{d}\bar{x}^{\frac{a}{2}}\|_{L^{2}}^{2};\nonumber\\
 K_{5}\leq & C \int |\nabla\mathbf{d}|^{3} |\nabla^{2}\mathbf{d}|\bar{x}^{a}\text{d}x
 \leq C\|\nabla^{2}\mathbf{d}\bar{x}^{\frac{a}{2}}\|_{L^{2}}
 \|\nabla\mathbf{d}\bar{x}^{\frac{a}{2}}\|_{L^{4}} \|\nabla \mathbf{d}\|_{L^{8}}^{2}\nonumber\\
 \leq &  C\|\nabla^{2}\mathbf{d}\bar{x}^{\frac{a}{2}}\|_{L^{2}}^2+
 C\|\nabla \mathbf{d}\bar{x}^{\frac{a}{2}}\|_{L^{2}}
 (\|\nabla^{2}\mathbf{d}\bar{x}^{\frac{a}{2}}\|_{L^{2}}+\|\nabla \mathbf{d}\nabla\bar{x}^{\frac{a}{2}}\|_{L^{2}}) \|\nabla \mathbf{d}\|_{L^{2}} \|\nabla^{2} \mathbf{d}\|_{L^{2}}^{3}\nonumber\\
 \leq& C(1+ \|\nabla^{2}\mathbf{d}\bar{x}^{\frac{a}{2}}\|_{L^{2}}^{2});\nonumber\\
 K_{6}\leq & C \int |\nabla\mathbf{d}||\nabla^{2} \mathbf{d}|\nabla\bar{x}^{a}\text{d}x\leq \|\nabla \mathbf{d}\bar{x}^{\frac{a}{2}}\|_{L^{2}}
 \|\nabla^{2}\mathbf{d}\bar{x}^{\frac{a}{2}}\|_{L^{4}}^{2}
 \nonumber\\
 \leq & C\|\nabla^{2}\mathbf{d}\bar{x}^{\frac{a}{2}}\|_{L^{2}}
 (\|\nabla^{3}\mathbf{d}\bar{x}^{\frac{a}{2}}\|_{L^{2}}
 +\|\nabla^{2}\mathbf{d}\bar{x}^{\frac{a}{2}}\|_{L^{2}})\nonumber\\
 \leq & \frac{1}{20}\|\nabla^{3}\mathbf{d}\bar{x}^{\frac{a}{2}}\|_{L^{2}}^{2}
 +C\|\nabla^{2}\mathbf{d}\bar{x}^{\frac{a}{2}}\|_{L^{2}}^{2};\nonumber\\
 K_{7}\leq & C\|\nabla \mathbf{u}\|_{L^{4}}^{2}\|\nabla\mathbf{d}\bar{x}^{\frac{a}{2}}\|_{L^{4}}^{2}
 \leq C \|\nabla\mathbf{u}\|_{L^{2}}\|\nabla^{2} \mathbf{u}\|_{L^{2}}\|\nabla\mathbf{d}\bar{x}^{\frac{a}{2}}\|_{L^{2}} \left(\|\nabla^{2} \mathbf{d}\bar{x}^{\frac{a}{2}}\|_{L^{2}} +\|\nabla\mathbf{d}\nabla \bar{x}^{\frac{a}{2}}\|_{L^{2}} \right)\nonumber\\
 \leq & C\|\nabla^{2}\mathbf{u}\|_{L^{2}}^{2}+C\|\nabla^{2} \mathbf{d}\bar{x}^{\frac{a}{2}}\|_{L^{2}}^2+C;\nonumber\\
 K_{8}+ K_{9}\leq & C\|\nabla\mathbf{d}\|_{L^{8}}^{4} \|\nabla^{2}\mathbf{d}\bar{x}^{\frac{a}{2}}\|_{L^{4}}^{2}+C\|\nabla \mathbf{d}\|_{L^{4}}^{2}\|\nabla^{2}\mathbf{d}\bar{x}^{\frac{a}{2}}\|_{L^{4}}^{2}
 \nonumber\\
 \leq & C(\|\nabla\mathbf{d}\|_{L^{2}}\|\nabla^{2}\mathbf{d}\|_{L^{2}}^{2}+\|\nabla \mathbf{d}\|_{L^{2}}\|\nabla^{2}\mathbf{d}\|_{L^{2}}) \|\nabla^{2} \mathbf{d}\bar{x}^{\frac{a}{2}}\|_{L^{2}}
 (\|\nabla^{3}\mathbf{d}\bar{x}^{\frac{a}{2}}\|_{L^{2}}\!
 +\!\|\nabla^{2}\mathbf{d}\nabla\bar{x}^{\frac{a}{2}}\|_{L^{2}})\nonumber\\
 \leq & \frac{1}{20}\|\nabla^{3}\mathbf{d}\bar{x}^{\frac{a}{2}}\|_{L^{2}}^{2}
 +C\|\nabla^{2}\mathbf{d}\bar{x}^{\frac{a}{2}}\|_{L^{2}}^{2};
 \nonumber\\
 K_{10}\leq & C\int |\nabla^{2}\mathbf{d}|^{2}\bar{x}^a \bar{x}^{-2}\log^{4}(e+|x|^2)\text{d}x
 \leq C\|\nabla^{2}\mathbf{d}\bar{x}^{\frac{a}{2}}\|_{L^{2}}^{2}.
 \end{align*}
Substituting the above estimates of
 $K_{i}\ (i=1,2,\cdots,10)$ into \eqref{eq3.37}, after by using \eqref{eq3.23}, we have
 \begin{align*}
 \frac{d}{dt} \|\nabla^{2}\mathbf{d}\bar{x}^{\frac{a}{2}}\|_{L^{2}}^{2}+&\! \|\nabla^{3} \mathbf{d}\bar{x}^{\frac{a}{2}}\|_{L^{2}}^{2}
  \leq C\|\nabla^{2}\mathbf{d}\bar{x}^{\frac{a}{2}}\|_{L^{2}}^{2}+C\|\nabla^{2} \mathbf{u}\|_{L^{2}}^{2}+C,
  \end{align*}
which multiplied by $t$ implies \eqref{eq3.35} after using  Gronwall's inequality, \eqref{eq3.26}, and \eqref{eq3.34}. This completes the proof of Lemma \ref{lem3.6}.  \hfill $\Box$

\begin{lemma}\label{lem3.7}
There exists a positive constant $C$ depending on $T$ such that
\begin{align}\label{eq3.38}
&\sup_{0\leq t\leq T} t\left(\|\sqrt{\rho} \mathbf{u}_{t}\|_{L^{2}}^{2}+\|\mathbf{d}_{t}\|_{H^{1}}^{2}+\|\nabla^{3} \mathbf{d}\|_{L^{2}}^{2}\right)+
\int_{0}^{T} t\left(\|\nabla^{2} \mathbf{u}_{t}\|_{L^{2}}^{2}+\|\nabla \mathbf{d}_{t}\|_{H^{1}}^{2}\right)\text{d}t
\leq C(T).
\end{align}
\end{lemma}

{\it Proof.}
First, we shall prove that
\begin{align}\label{eq3.39}
\int_{0}^{T} \left(\|\sqrt{\rho} \mathbf{u}_{t}\|_{L^{2}}^{2}+\|\nabla \mathbf{d}_{t}\|_{L^{2}}^{2}\right)\text{d}t\leq C(T).
\end{align}
On the one hand, we derive from H{\"o}lder's inequality, \eqref{l2}, \eqref{eq3.2}, and \eqref{ll} that
\begin{align}\label{eq3.40}
\|\sqrt{\rho} \mathbf{u}_{t}\|_{L^{2}}^{2}
\leq & \|\sqrt{\rho}\dot{\mathbf{u}}\|_{L^{2}}^{2}+\|\sqrt{\rho}|\mathbf{u}||\nabla \mathbf{u}|\|_{L^{2}}^{2}\nonumber\\
\leq & \|\sqrt{\rho}\dot{\mathbf{u}}\|_{L^{2}}^{2}+C\|\sqrt{\rho} \mathbf{u}\|_{L^{6}}^{2} \|\nabla \mathbf{u}\|_{L^{3}}^{2}\nonumber\\
\leq &  \|\sqrt{\rho}\dot{\mathbf{u}}\|_{L^{2}}^{2}+C\left(\|\sqrt{\rho} \mathbf{u}\|_{L^{2}}^{2}+\|\nabla \mathbf{u}\|_{L^{2}}^{2}\right) \|\nabla \mathbf{u}\|_{L^{2}}^{\frac{4}{3}} \|\nabla^{2} \mathbf{u}\|_{L^{2}}^{\frac{2}{3}}\nonumber\\
\leq & \|\sqrt{\rho}\dot{\mathbf{u}}\|_{L^{2}}^{2}+C\left(1+\|\nabla^{2} \mathbf{u}\|_{L^{2}}^{2}\right).
\end{align}
On the other hand, by virtue of \eqref{eq3.7}, \eqref{l3}, \eqref{eq3.2}, \eqref{ll}, and \eqref{eq3.34}, we obtain
\begin{align}\label{eq3.41}
\|\nabla \mathbf{d}_{t}\|_{L^{2}}^{2}
\leq & C\left(\|\nabla^{3} \mathbf{d}\|_{L^{2}}^{2}+\||\nabla\mathbf{u}| |\nabla \mathbf{d}|\|_{L^{2}}^{2}+\||\mathbf{u}||\nabla^{2}\mathbf{d}|\|
_{L^{2}}^{2}+\||\nabla \mathbf{d}|^{3}\|_{L^{2}}^{2}+\||\nabla \mathbf{d}||\nabla^{2}\mathbf{d}|\|_{L^{2}}^{2}\right)\nonumber\\
\leq & C\left(\|\nabla^{3}\mathbf{d}\|_{L^{2}}^{2}+\|\nabla\mathbf{u}\|_{L^{4}}^{2} \|\nabla \mathbf{d}\|_{L^{4}}^{2}+\|\mathbf{u}\bar{x}^{-\frac{a}{4}}\|_{L^{8}}^{2}
\|\nabla^{2}\mathbf{d}\bar{x}^{\frac{a}{2}}\|_{L^{2}}
\|\nabla^{2}\mathbf{d}\|_{L^{4}}\right)\nonumber\\
& +C\left(\|\nabla\mathbf{d}\|_{L^{6}}^{6}+\||\nabla \mathbf{d}||\nabla^{2}\mathbf{d}|\|_{L^{2}}^{2}\right)\nonumber\\
\leq & C\left(\|\nabla^{3}\mathbf{d}\|_{L^{2}}^{2}+ \|\nabla^{2} \mathbf{u}\|_{L^{2}}^{2}+\|\nabla^{2}\mathbf{d}\bar{x}^{\frac{a}{2}}\|_{L^{2}}^{2}
+\|\mathbf{u}\bar{x}^{-\frac{a}{4}}\|_{L^{8}}^{4}
\|\nabla^{2}\mathbf{d}\|_{L^{4}}^{2}+1\right)\nonumber\\
\leq & C\left(\|\nabla^{3} \mathbf{d}\|_{L^{2}}^{2}+ \|\nabla^{2} \mathbf{u}\|_{L^{2}}^{2}+\|\nabla^{2}\mathbf{d}\bar{x}^{\frac{a}{2}}\|_{L^{2}}^{2}
+(\|\sqrt{\rho}\mathbf{u}\|_{L^{2}}^{2}+\|\nabla\mathbf{u}\|_{L^{2}}^{2})^{2}
\|\nabla^{2}\mathbf{d}\|_{L^{4}}^{2}+1\right)\nonumber\\
\leq & C\left(\|\nabla^{3} \mathbf{d}\|_{L^{2}}^{2}+ \|\nabla^{2}\mathbf{u}\|_{L^{2}}^{2}+1\right),
\end{align}
which combined with \eqref{eq3.40}, \eqref{eq3.2}, and \eqref{eq3.26} leads to \eqref{eq3.39}.

Now, differentiating $\eqref{eq1.1}_{2}$ with respect to $t$ gives
\begin{align*}
\rho \mathbf{u}_{tt} +\rho \mathbf{u}\cdot\nabla \mathbf{u}_{t}-\Delta \mathbf{u}_{t}+\nabla P_{t}= -\rho_{t}(\mathbf{u}_{t}+\mathbf{u}\cdot\nabla\mathbf{u})-\rho \mathbf{u}_{t}\cdot\nabla \mathbf{u} -\divv(\nabla \mathbf{d}\odot \nabla \mathbf{d})_{t}.
\end{align*}
Multiplying the above equality by $\mathbf{u}_{t}$ and integrating the resulting equality by parts over $\mathbb{R}^{2}$, we deduce after using $\eqref{eq1.1}_{1}$ and $\eqref{eq1.1}_{4}$ that
\begin{align}\label{eq3.42}
\frac{1}{2}\frac{d}{dt} \int \rho|\mathbf{u}_t|^2dx
+\int|\nabla \mathbf{u}_t|^2dx
\leq & C\int\! \rho |\mathbf{u}| |\mathbf{u}_{t}|\left( |\nabla\mathbf{u}_{t}|\!+ \!|\nabla \mathbf{u}|^{2}\! +\!|\mathbf{u}||\nabla^{2}\mathbf{u}|\right)\text{d}x+C\int \rho |\mathbf{u}|^{2}|\nabla\mathbf{u}| |\nabla\mathbf{u}_{t}|\text{d}x
\nonumber\\
& + C\int \rho |\mathbf{u}_{t}|^{2}|\nabla\mathbf{u}|\text{d}x +C\int|\nabla\mathbf{d}||\nabla\mathbf{d}_{t}|
|\nabla\mathbf{u}_{t}|\text{d}x\nonumber\\
\triangleq & M_{1}+M_{2}+M_{3}+M_{4}.
\end{align}
The terms on the right-hand side of \eqref{eq3.42} can be bounded as follows.

By \eqref{l3}, \eqref{l2}, \eqref{eq3.2}, \eqref{ll}, H\"{o}lder's and Gagliardo-Nirenberg inequalities, we have
\begin{align*}
M_{1}& \leq C\|\sqrt{\rho}\mathbf{u}\|_{L^{6}} \|\sqrt{\rho} \mathbf{u}_{t}\|_{L^{2}}^{\frac{1}{2}}
\|\sqrt{\rho} \mathbf{u}_{t}\|_{L^{6}}^{\frac{1}{2}} \left(\|\nabla \mathbf{u}_{t}\|_{L^{2}}+\|\nabla\mathbf{u}\|_{L^{4}}^{2}\right)\\
& \quad +C\|\rho^{\frac{1}{4}}\mathbf{u}\|_{L^{12}}^{2} \|\sqrt{\rho} \mathbf{u}_{t}\|_{L^{2}}^{\frac{1}{2}}
\|\sqrt{\rho}\mathbf{u}_{t}\|_{L^{6}}^{\frac{1}{2}} \|\nabla^{2} \mathbf{u}\|_{L^{2}}\nonumber\\
& \leq C\|\sqrt{\rho} \mathbf{u}_{t}\|_{L^{2}}^{\frac{1}{2}}
\left(\|\sqrt{\rho}\mathbf{u}_{t}\|_{L^{2}}\! +\!\|\nabla \mathbf{u}_{t}\|_{L^{2}}\right)^{\frac{1}{2}} \left(\|\nabla\mathbf{u}_{t}\|_{L^{2}}
+\!\|\nabla^{2}\mathbf{u}\|_{L^{2}}\right)\nonumber\\
& \leq \frac{1}{6} \|\nabla\mathbf{u}_{t}\|_{L^{2}}^{2}+ C\left(1+\|\sqrt{\rho}\mathbf{u}_{t}\|_{L^{2}}^{2}
+\|\nabla^{2}\mathbf{u}\|_{L^{2}}^{2}\right).
\end{align*}
Next, H{\"o}lder's inequality, \eqref{l3}, and \eqref{l2} imply
\begin{align*}
M_{2}+M_{3} \leq & C \|\sqrt{\rho}\mathbf{u}\|_{L^{8}}^{2} \|\nabla \mathbf{u}\|_{L^{4}}\|\nabla\mathbf{u}_{t}\|_{L^{4}}+C\|\sqrt{\rho} \mathbf{u}_{t}\|_{L^{6}}^{\frac{3}{2}} \|\sqrt{\rho} \mathbf{u}_{t}\|_{L^{2}}^{\frac{1}{2}}\|\nabla\mathbf{u}\|_{L^{2}}\nonumber\\
\leq & \frac{1}{6} \|\nabla\mathbf{u}_{t}\|_{L^{2}}^{2}+ C\left(1+\|\sqrt{\rho} \mathbf{u}_{t}\|_{L^{2}}^{2}
+\|\nabla^{2}\mathbf{u}\|_{L^{2}}^{2}\right).
\end{align*}
For the term $M_4$, by Ladyzhenskaya's inequality, \eqref{eq3.2}, and \eqref{ll},  we have
\begin{align*}
M_{4}\leq & C\|\nabla\mathbf{d}\|_{L^{4}} \|\nabla\mathbf{d}_{t}\|_{L^{4}} \|\nabla\mathbf{u}_{t}\|_{L^{2}}\\
\leq & \frac{1}{6} \|\nabla\mathbf{u}_{t}\|_{L^{2}}^{2} +C\|\nabla\mathbf{d}\|_{L^{2}} \|\nabla^2\mathbf{d}\|_{L^{2}} \|\nabla\mathbf{d}_{t}\|_{L^{2}}\|\nabla^{2} \mathbf{d}_{t}\|_{L^{2}}
\\
\leq &\frac{1}{6} \|\nabla \mathbf{u}_{t}\|_{L^{2}}^{2} +\frac{1}{4(C_{3}+1)}\|\nabla^{2} \mathbf{d}_{t}\|_{L^{2}}^{2}+C\|\nabla\mathbf{d}_{t}\|_{L^{2}}^{2},
\end{align*}
where the positive constant $C_{3}$ is defined in the following \eqref{eq3.44} and \eqref{eq3.47}. Substituting the estimates of $M_{i}\ (i=1,2,\cdots, 4)$ into \eqref{eq3.42},  there holds
\begin{align}\label{eq3.43}
\frac{d}{dt} \|\sqrt{\rho} \mathbf{u}_{t}\|_{L^{2}}^{2}+\|\nabla \mathbf{u}_{t}\|_{L^{2}}^{2}
\leq C(\|\sqrt{\rho}\mathbf{u}_{t}\|_{L^{2}}^{2}
+ \|\nabla \mathbf{d}_{t}\|_{L^{2}}^{2})
+ \frac{1}{2(C_{3}+1)}\|\nabla^{2} \mathbf{d}_{t}\|_{L^{2}}^{2}+C(\|\nabla^{2}\mathbf{u}\|_{L^{2}}^{2}+1).
\end{align}

Next, differentiating $\eqref{eq1.1}_{3}$ with respect to $t$, and multiplying the resulting equality with $\mathbf{d}_{t}$ and then integrating by parts over $\mathbb{R}^{2}$, we arrive at
\begin{align*}
\frac{1}{2}\frac{d}{dt} \int|\mathbf{d}_{t}|^2 \text{d}x
+ \int|\nabla \mathbf{d}_{t}|^2 \text{d}x
\leq & C\int |\mathbf{u}_{t}||\nabla \mathbf{d}||\mathbf{d}_{t}|\text{d}x +C\int |\nabla \mathbf{d}_{t}| |\nabla \mathbf{d}||\mathbf{d}_{t}|\text{d}x+C\int |\nabla\mathbf{d}|^{2} |\mathbf{d}_{t}|^{2}\text{d}x \\
\triangleq & M_{5}+M_{6}+M_{7}.
\end{align*}
By H\"{o}lder's and Ladyzhenskaya's inequalities, \eqref{l3}, \eqref{eq3.2}, \eqref{ll}, and \eqref{eq3.34}, we derive
\begin{align*}
M_{5}\leq & C\|\mathbf{u}_{t}\bar{x}^{-\frac{a}{2}}\|_{L^{4}} \|\nabla\mathbf{d} \bar{x}^{\frac{a}{2}}\|_{L^{2}} \|\mathbf{d}_{t}\|_{L^{4}}
\\
\leq & C(\|\sqrt{\rho} \mathbf{u}_{t}\|_{L^{2}}+\|\nabla\mathbf{u}_{t}\|_{L^{2}})\|\nabla \mathbf{d}\bar{x}^{\frac{a}{2}}\|_{L^{2}} \|\mathbf{d}_{t}\|_{L^{2}}^{\frac{1}{2}}\| \nabla \mathbf{d}_{t}\|_{L^{2}}^{\frac{1}{2}}\nonumber\\
\leq & \frac{1}{4}\|\nabla \mathbf{u}_{t}\|_{L^{2}}^{2}+C(\|\sqrt{\rho} \mathbf{u}_{t}\|_{L^{2}}^{2}+ \|\mathbf{d}_{t}\|_{L^{2}}^{2}+\|\nabla \mathbf{d}_{t}\|_{L^{2}}^{2});\nonumber\\
M_{6}+M_{7}\leq &\frac{1}{4}\|\nabla \mathbf{d}_{t}\|_{L^{2}}+C\|\nabla \mathbf{d}\|_{L^{4}}^{2}\|\mathbf{d}_{t}\|_{L^{4}}^{2}
\nonumber\\
\leq &  \frac{1}{4}\|\nabla \mathbf{d}_{t}\|_{L^{2}}+C\|\nabla \mathbf{d}\|_{L^{2}}\|\nabla^{2} \mathbf{d}\|_{L^{2}}\|\mathbf{d}_{t}\|_{L^{2}}\| \nabla \mathbf{d}_{t}\|_{L^{2}}\nonumber\\
\leq & \frac{1}{2}\|\nabla \mathbf{d}_{t}\|_{L^{2}}+C\|\mathbf{d}_{t}\|_{L^{2}}^{2}.
\end{align*}
Hence
\begin{align}\label{eq3.44}
\frac{d}{dt} \|\mathbf{d}_{t}\|_{L^{2}}^{2}+\|\nabla\mathbf{d}_{t}\|_{L^{2}}^{2}\leq& C_{3}\|\nabla \mathbf{u}_{t}\|_{L^{2}}^{2}+C\left(\|\sqrt{\rho}\mathbf{u}_{t}\|_{L^{2}}^{2}+ \|\mathbf{d}_{t}\|_{L^{2}}^{2}+\|\nabla \mathbf{d}_{t}\|_{L^{2}}^{2}\right).
\end{align}

Differentiating \eqref{eq3.7} with respect to time variable $t$ ensures
\begin{align}\label{eq3.45}
\nabla \mathbf{d}_{tt}-\Delta\nabla \mathbf{d}_{t} =-\nabla (\mathbf{u}\cdot \nabla \mathbf{d})_{t}+\nabla (|\nabla \mathbf{d}|^{2}\mathbf{d})_{t}.
\end{align}
Multiplying \eqref{eq3.45} by $\nabla \mathbf{d}_{t}$, and integrating the resulting equality over $\mathbb{R}^2$, we find
\begin{align}\label{eq3.46}
\frac{1}{2}\frac{d}{dt}\|\nabla \mathbf{d}_{t}\|_{L^{2}}^{2}+\|\nabla^{2} \mathbf{d}_{t}\|_{L^{2}}^{2}
\leq & C\int |\nabla \mathbf{u}_{t}||\nabla \mathbf{d}||\nabla \mathbf{d}_{t}|\text{d}x+C\int |\nabla \mathbf{u}||\nabla \mathbf{d}_{t}|^{2}\text{d}x
+C\int |\mathbf{u}_{t}||\nabla^{2} \mathbf{d}| |\nabla \mathbf{d}_{t}|\text{d}x\nonumber\\
& +C\int |\nabla \mathbf{d}|^{2}|\mathbf{d}_{t}||\nabla^{2} \mathbf{d}_{t}|\text{d}x+C\int |\nabla \mathbf{d}||\nabla \mathbf{d}_{t}| |\nabla^{2} \mathbf{d}_{t}|\text{d}x\nonumber\\
\triangleq & M_{8}+M_{9}+M_{10}+M_{11}+M_{12}.
\end{align}
By H\"{o}lder's and Ladyzhenskaya's inequalities, \eqref{eq3.2}, and \eqref{ll}, we have
\begin{align*}
M_{8}\leq & C\|\nabla \mathbf{u}_{t}\|_{L^{2}}\|\nabla \mathbf{d}_{t}\|_{L^{4}}\|\nabla \mathbf{d}\|_{L^{4}}\nonumber\\
\leq & \frac{1}{2} \|\nabla \mathbf{u}_{t}\|_{L^{2}}^{2}+C
 \|\nabla \mathbf{d}_{t}\|_{L^{2}}\|\nabla^{2} \mathbf{d}_{t}\|_{L^{2}}\|\nabla \mathbf{d}\|_{L^{2}}\|\nabla^{2} \mathbf{d}\|_{L^{2}}\nonumber\\
\leq & \frac{1}{2} \|\nabla \mathbf{u}_{t}\|_{L^{2}}^{2}+ \frac{1}{4} \|\nabla^{2} \mathbf{d}_{t} \|_{L^{2}}^{2}
+C\|\nabla \mathbf{d}_{t}\|_{L^{2}}^{2}.
\end{align*}
Similarly, we get
\begin{align*}
M_{9}\leq & C\|\nabla \mathbf{u}\|_{L^{2}}\|\nabla \mathbf{d}_{t}\|_{L^{4}}^{2}
\leq C\|\nabla \mathbf{u}\|_{L^{2}}\|\nabla \mathbf{d}_{t}\|_{L^{2}} \|\nabla^{2} \mathbf{d}_{t}\|_{L^{2}}
\leq \frac{1}{16} \|\nabla^{2} \mathbf{d}_{t}\|_{L^{2}}^2 +C\|\nabla \mathbf{d}_{t}\|_{L^{2}}^2;\\
M_{12}\leq & \frac{1}{32}\|\nabla^{2} \mathbf{d}_{t}\|_{L^{2}}^{2}+C\|\nabla \mathbf{d}\|_{L^{4}}^{2}\| \nabla \mathbf{d}_{t}\|_{L^{4}}^{2}\leq \frac{1}{16}\|\nabla^{2} \mathbf{d}_{t}\|_{L^{2}}^{2}+C\|\nabla \mathbf{d}_{t}\|_{L^{2}}^{2}.
\end{align*}
Applying H\"{o}lder's and Ladyzhenskaya's inequalities, \eqref{l3}, \eqref{eq3.2}, \eqref{ll}, and \eqref{eq3.34}, we obtain
\begin{align*}
M_{10}\leq &
C\|\mathbf{u}_{t}\bar{x}^{-\frac{a}{4}}\|_{L^{4}}
\|\nabla^{2}\mathbf{d}\bar{x}^{\frac{a}{2}}\|_{L^{2}}^{\frac{1}{2}} \|\nabla^{2}\mathbf{d}\|_{L^{2}}^{\frac{1}{2}}
\|\nabla \mathbf{d}_{t}\|_{L^{4}}\nonumber\\
\leq &  C\|\mathbf{u}_{t}\bar{x}^{-\frac{a}{4}}\|_{L^{4}}\|\nabla^{2}\mathbf{d}
\bar{x}^{\frac{a}{2}}\|_{L^{2}}^{\frac{1}{2}} \|\nabla^{2}\mathbf{d}\|_{L^{2}}^{\frac{1}{2}}
\|\nabla \mathbf{d}_{t}\|_{L^{2}}^{\frac{1}{2}}\|\nabla^{2} \mathbf{d}_{t}\|_{L^{2}}^{\frac{1}{2}}\nonumber\\
\leq & \frac{1}{16} \|\nabla^{2}\mathbf{d}_{t}\|_{L^{2}}^{2}+C\|\mathbf{u}_{t}\bar{x}
^{-\frac{a}{4}}\|_{L^{4}}^{2}
+C_{1}(t)\|\nabla \mathbf{d}_{t}\|_{L^{2}}^{2}\nonumber\\
\leq & \frac{1}{16} \|\nabla^{2}\mathbf{d}_{t}\|_{L^{2}}^{2}
+C(\|\sqrt{\rho} \mathbf{u}_{t}\|_{L^{2}}^2+\|\nabla \mathbf{u}_{t}\|_{L^{2}}^{2})+C_{1}(t)
\|\nabla \mathbf{d}_{t}\|_{L^{2}}^{2},
\end{align*}
where $C_{1}(t)\geq0,\ \int_{0}^{T}C_{1}(t)dt\leq C(T)$ (for all $T\in(0,\infty)$).
By H\"{o}lder's and Gagliardo-Nirenberg inequalities, \eqref{eq3.2}, \eqref{ll}, and \eqref{eq3.34}, we deduce
\begin{align*}
M_{11}\leq & \frac{1}{16}\|\nabla^{2} \mathbf{d}_{t}\|_{L^{2}}^{2}+C\|\nabla \mathbf{d}\|_{L^{8}}^{2}\| \mathbf{d}_{t}\|_{L^{4}}^{2}
\nonumber\\
\leq &  \frac{1}{16}\|\nabla^{2} \mathbf{d}_{t}\|_{L^{2}}^{2}+C\|\nabla \mathbf{d}\|_{L^{2}}^{\frac{1}{2}}\|\nabla^{2} \mathbf{d}\|_{L^{2}}^{\frac{3}{2}}\| \mathbf{d}_{t}\|_{L^{2}}\|\nabla \mathbf{d}_{t}\|_{L^{2}}\nonumber\\
\leq &  \frac{1}{16}\|\nabla^{2} \mathbf{d}_{t}\|_{L^{2}}^{2}+C(\| \mathbf{d}_{t}\|_{L^{2}}^{2}+\|\nabla \mathbf{d}_{t}\|_{L^{2}}^{2}).
\end{align*}
Inserting the estimates of $M_{i}\ (i=8,9,\cdots, 12)$ into \eqref{eq3.46}, it follows that
\begin{align}\label{eq3.47}
&\frac{d}{dt}\|\nabla \mathbf{d}_{t}\|_{L^{2}}^{2}+\|\nabla^{2} \mathbf{d}_{t}\|_{L^{2}}^{2}
\leq C_{3} \left(\|\sqrt{\rho} \mathbf{u}_{t}\|_{L^{2}}+\|\nabla \mathbf{u}_{t}\|_{L^{2}}^{2}\right)+C
 \left(\|\mathbf{d}_{t}\|_{L^{2}}^{2}+\|\nabla \mathbf{d}_{t}\|_{L^{2}}^{2}\right).
\end{align}

Now, multiplying \eqref{eq3.43} by  ${2(C_{3}+1)}$ and adding the resulting inequality with \eqref{eq3.44} and  \eqref{eq3.47}, we infer that
\begin{align*}
&\frac{d}{dt}({2(C_{3}+1)}\|\sqrt{\rho} \mathbf{u}_{t}\|_{L^{2}}^{2}+\|\mathbf{d}_{t}\|_{H^{1}}^{2})+\|\nabla \mathbf{u}_{t}\|_{L^{2}}
+\frac{1}{2} \|\nabla \mathbf{d}_{t}\|_{H^{1}}^{2} \nonumber\\
& \leq C\left(1+\|\sqrt{\rho} \mathbf{u}_{t}\|_{L^{2}}^{2}+\|\nabla \mathbf{d}_{t}\|_{H^{1}}^{2}\right)
+C\left(1+\|\nabla^{2} \mathbf{u}\|_{L^{2}}^{2}\right),
\end{align*}
which multiplied  by $t$, together with Gronwall's inequality, \eqref{eq3.39}, and \eqref{eq3.10} yields
\begin{align}\label{eq3.48}
\sup_{0\leq t\leq T} t\left(\|\sqrt{\rho} \mathbf{u}_{t}\|_{L^{2}}^{2}+\|\mathbf{d}_{t}\|_{H^{1}}^{2}\right)
+\int_{0}^{T} t\left(\|\nabla \mathbf{u}_{t}\|_{L^{2}}^{2}+\|\nabla \mathbf{d}_{t}\|_{H^{1}}^{2}\right)\text{d}t\leq C(T).
\end{align}

Finally, it follows from \eqref{eq3.7}, H\"{o}lder's and Gagliardo-Nirenberg inequalities, \eqref{l3}, \eqref{eq3.2}, \eqref{ll}, and $|\mathbf{d}|=1$  that
\begin{align*}
\|\nabla^{3}\mathbf{d}\|_{L^{2}}^{2}
\leq & C\left(\|\nabla \mathbf{d}_{t}\|_{L^{2}}^{2}+\||\nabla \mathbf{u}||\nabla \mathbf{d}|\|_{L^{2}}^{2}
+\||\mathbf{u}||\nabla^{2}\mathbf{d}|\|_{L^{2}}^{2}+\||\nabla \mathbf{d}|^{3}\|_{L^{2}}^{2}+\||\nabla^{2}\mathbf{d}||\nabla \mathbf{d}|\|_{L^{2}}^{2}\right)\nonumber\\
\leq & C(\|\nabla \mathbf{d}_{t}\|_{L^{2}}^{2}+\|\nabla \mathbf{u}\|_{L^{4}}^{2} \|\nabla \mathbf{d}\|_{L^{4}}^{2}+\|\mathbf{u}\bar{x}^{-\frac{a}{4}}\|_{L^{8}}^{2}
\|\nabla^{2}\mathbf{d}\bar{x}^{\frac{a}{2}}\|_{L^{2}}
\|\nabla^{2}\mathbf{d}\|_{L^{4}})\nonumber\\
& + C(\|\nabla \mathbf{d}\|_{L^{6}}^{6}+\||\nabla \mathbf{d}||\nabla^{2}\mathbf{d}|\|_{L^{2}}^{2})\nonumber\\
\leq & C(\|\nabla \mathbf{d}_{t}\|_{L^{2}}^{2}+ \|\nabla^{2} \mathbf{u}\|_{L^{2}}^{2}+\|\nabla^{2}\mathbf{d}\bar{x}^{\frac{a}{2}}\|_{L^{2}}^{2}
+\|\mathbf{u}\bar{x}^{-\frac{a}{4}}\|_{L^{8}}^{4}
\|\nabla^{2}\mathbf{d}\|_{L^{4}}^{2}\nonumber\\
& +\|\nabla \mathbf{d}\|_{L^{2}}^{2}\|\nabla^{2}\mathbf{d}\|_{L^{2}}^{4}
+\|\nabla^{2} \mathbf{d}\|_{L^{3}}^{2}\|\nabla \mathbf{d}\|_{L^{6}}^{2})\nonumber\\
\leq & C(\|\nabla \mathbf{d}_{t}\|_{L^{2}}^{2}\!+ \!\|\nabla^{2} \mathbf{u}\|_{L^{2}}^{2}\!+\!\|\nabla^{2}\mathbf{d}\bar{x}^{\frac{a}{2}}\|_{L^{2}}^{2}\!
+\!(\|\sqrt{\rho} \mathbf{u}\|_{L^{2}}^{2}\!+\!\|\nabla \mathbf{u}\|_{L^{2}}^{2})^{2}
\|\nabla^{2}\mathbf{d}\|_{L^{2}}\|\nabla^{3}\mathbf{d}\|_{L^{2}}\nonumber\\
&+\|\nabla^{2}\mathbf{d}\|_{L^{2}}^{4}+\|\nabla^{2}\mathbf{d}\|_{L^{2}}
\|\nabla^{3}\mathbf{d}\|_{L^{2}})\nonumber\\
\leq &\frac{1}{2}\|\nabla^{3}\mathbf{d}\|_{L^{2}}^{2}
+C\left(\|\nabla \mathbf{d}_{t}\|_{L^{2}}^{2}+ \|\nabla^{2} \mathbf{u}\|_{L^{2}}^{2}+\|\nabla^{2}\mathbf{d}\bar{x}^{\frac{a}{2}}\|_{L^{2}}^{2}
+1\right),
\end{align*}
which combined with \eqref{eq3.3}, \eqref{eq3.14}, \eqref{eq3.35}, and \eqref{eq3.48} implies \eqref{eq3.38}.
The proof of Lemma \ref{lem3.7} is completed. \hfill$\Box$

\section{Proof of Theorem \ref{thm1.2}}
With the a priori estimates in Section \ref{sec3} in hand, we are now in a position to prove Theorem \ref{thm1.2}.

By Lemma \ref{lem2.1}, we know that there exists a $T_{*}>0$ such that the Cauchy problem of system \eqref{eq1.1}--\eqref{eq1.2} admits a unique strong solution $(\rho, \mathbf{u},P, \mathbf{d})$ on $\mathbb{R}^{2}\times (0,T_{*}]$. In what follows, we shall extend the local solution to all the time.

Set
\begin{align}\label{eq4.1}
T^{*}\!=\sup \left\{ T | (\rho,\mathbf{u}, P, \mathbf{d}) \text{ is a strong solution to \eqref{eq1.1}--\eqref{eq1.2} on } \mathbb{R}^{2}\times (0,T] \right\}.
\end{align}
First, for any $0<\tau<T_{*}<T\leq T^{*}$ with $T$ finite, one deduces from \eqref{eq3.2}, \eqref{ll}, \eqref{eq3.14}, and \eqref{eq3.38} that for all $q\geq 2$,
\begin{align}\label{eq4.2}
\nabla \mathbf{u}, \nabla \mathbf{d},\nabla^{2}\mathbf{d}\in C([\tau,T]; L^{2}\cap L^{q}),
\end{align}
where one has used the standard embedding
\begin{align*}
L^{\infty}(\tau,T; H^{1})\cap H^{1}(\tau,T; H^{-1}) \hookrightarrow C(\tau,T;L^{q})\quad \text{ for all }q\in [2,\infty).
\end{align*}
Moreover, it follows from \eqref{eq3.26}, \eqref{eq3.33}, and \cite[Lemma 2.3]{Lions1} that
\begin{align}\label{eq4.3}
\rho\in C([0,T];L^{1}\cap H^{1}\cap W^{1,q}).
\end{align}

Now, we claim that
\begin{align}\label{eq4.4}
T^{*}=\infty.
\end{align}
Otherwise, if $T^{*}<\infty$, it follows from \eqref{eq4.2}, \eqref{eq4.3}, \eqref{eq3.2}, \eqref{ll}, \eqref{eq3.33}, and \eqref{eq3.34} that
\begin{align*}
(\rho,\mathbf{u},P,\mathbf{d})(x,T^{*})= \lim_{t\rightarrow T^{*}}(\rho,\mathbf{u},P,\mathbf{d})(x,t)
\end{align*}
satisfies the initial conditions \eqref{eq1.6} at $t=T^{*}$. Moreover, by using \eqref{eq1.4} and \eqref{eq3.5} with $p=1$, it follows that
\begin{align}
\int\rho(x,T^{*})\text{d}x=\int \rho_{0}(x) \text{d}x=1,
\end{align}
and  notice that there exists  $N_{0}>0$,
it is easy to see that
\begin{align*}
\int \rho(x,T^{*})\text{d}x\geq \frac{1}{2} \int_{B_{N_{0}}}\rho(x,T^{*})\text{d}x\geq \frac{1}{2}.
\end{align*}
Thus, we can take $(\rho,\mathbf{u},P,\mathbf{d})(x,T^{*})$ as the initial data, Lemma \ref{lem2.1} implies that one could extend the local solutions beyond $T^{*}$. This contradicts the assumption of $T^{*}$ in \eqref{eq4.1}. Hence, we prove \eqref{eq4.4}. Furthermore, from \eqref{eq3.3}, \eqref{eq3.13}, and \eqref{eq3.14}, one obtains that \eqref{eq1.9} holds.
This completes the proof of Theorem \ref{thm1.2}.
\hfill$\Box$

\section*{Acknowledgments}
The authors thank Dr. Jinkai Li for introducing this topic and for all the helpful discussions. Part of the work was done when Xin Zhong was visiting The Institute of Mathematical Sciences, The Chinese University of Hong
Kong. He thanks Professor Zhouping Xin for the invitation and constant help during the visit. Lin Li is supported by National Natural Science Foundation of China (11601046), Chongqing Science and Technology Commission (cstc2016jcyjA0310), Chongqing Municipal Education Commission (KJ1600603) and Program for University Innovation Team of Chongqing (CXTDX201601026). Qiao Liu is partially supported by National Natural Science Foundation of China (11401202) and China Postdoctoral Science Foundation
(2015M570053, 2016T90063). Xin Zhong is supported by Fundamental Research Funds for the Central Universities (XDJK2017C050),
China Postdoctoral Science Foundation (2017M610579), and the Doctoral Fund of Southwest University (SWU116033). At the same time, the authors express their gratitude to the reviewers for careful reading and helpful suggestions which led to an improvement of the original manuscript.

\end{document}